\documentclass[12pt,reqno]{amsart}   	
\usepackage[letterpaper, margin=1in]{geometry}
\usepackage{graphicx}				
\usepackage{amssymb}
\usepackage{mathtools,amsmath}  
\usepackage{amsthm}
\usepackage{amssymb}
\usepackage{mathrsfs}
\usepackage{epstopdf}
\usepackage{color}
\usepackage{tikz-cd}
\usepackage{tikz}
\usepackage{tikz-3dplot}
\usepackage{subfigure}
\usepackage{enumerate}
\usepackage{comment}
\usepackage{braket}
\usepackage{algpseudocode}
\usepackage[pagebackref,colorlinks,citecolor=blue,linkcolor=magenta]{hyperref}

\usepackage[linesnumbered,ruled]{algorithm2e}
\RequirePackage{amsthm,amsmath,amsfonts,amssymb}
\usepackage[utf8]{inputenc}
\usepackage{chngcntr}
\usepackage{float}

\theoremstyle{plain}
   \newtheorem{theorem}{Theorem}[section]
   \newtheorem{proposition}[theorem]{Proposition}
   \newtheorem{lemma}[theorem]{Lemma}
   \newtheorem{corollary}[theorem]{Corollary}

\theoremstyle{definition}
   \newtheorem{definition}[theorem]{Definition}

   \newtheorem{example}[theorem]{Example}
   
\theoremstyle{remark}
 \newtheorem{remark}[theorem]{Remark}

\newcommand{\LL}{\mathcal{L}}

\newcommand{\HH}{\mathcal{H}}

\newcommand{\MM}{\mathcal{M}}

\newcommand{\Z}{\mathbb{Z}}

\newcommand{\R}{\mathbb{R}}
\newcommand{\C}{\mathbb{C}}

\newcommand{\pp}{\mathbb{P}}

\newcommand{\vol}{\operatorname{Vol}}

\newcommand{\mldeg}{\operatorname{mldeg}}
\newcommand{\pr}{\mathrm{pr}}
\newcommand{\ord}{\mathrm{ord}}

\SetKwInput{KwInput}{Input}
\SetKwInput{KwOutput}{Output}

\def\newop#1{\expandafter\def\csname #1\endcsname{\mathop{\rm
#1}\nolimits}}

\usetikzlibrary{arrows.meta,fit}
\newcommand\addvmargin{%
  \node[fit=(current bounding box),inner ysep=3pt,inner xsep=0pt] {};
}
\xdef\baselinedim{\number\dimexpr\baselineskip\relax sp}
\newcommand\addbase{%
  \node (base) {\phantom{1}};
}

\newlength{\figureheight}
\newsavebox{\figurecaptionbox}

\newop{Inv}
\newop{conv}
\newop{pa}
\newop{de}
\newop{an}
\newop{nd}
\newop{ch}
\newop{CS}

\begin{document}
\title[ML degrees of 3D polytopes built from polygons of ML degree one]{Maximum likelihood degrees of 3D polytopes built from polygons of ML degree one}
\author{Eliana Duarte}
\address{Centro de Matem\'atica da Universidade do Porto (CMUP), Departamento de Matem\'atica, Faculdade de Ci\^encias, Universidade do Porto, Rua do Campo Alegre s/n, 4169-007 Porto, Portugal}
\email{eliana.gelvez@fc.up.pt}

\author{Federico Lazzeri}
\address{Department of Mathematics, KTH Royal Institute of Technology, SE-100 44 Stockholm, Sweden}
\email{lazzeri@kth.se}

\author{Bowen Li}
\address{Department of Mathematics, University of Illinois Urbana-Champaign, 1409 West Green Street, Urbana, IL 61801, USA}
\email{bowenl6@illinois.edu}

\author{Ludovico Piazza}
\address{Dipartimento di Matematica, Universit\`a di Pisa, Largo Bruno Pontecorvo 5, 56127 Pisa, Italy}
\email{l.piazza5@studenti.unipi.it}

\author{Bernhard Reinke}
\address{Aix-Marseille Universit\'e, CNRS, Institut de Math\'ematiques de Marseille (I2M), Marseille, France}
\email{bernhard.reinke@univ-amu.fr}

\author{Nico Wolf}
\address{HTW Berlin -- University of Applied Sciences, Berlin, Germany}
\email{nico.wolf@htw-berlin.de}
\renewcommand{\shortauthors}{E. Duarte, F. Lazzeri, B. Li, L. Piazza, B. Reinke, N. Wolf}
\date{\today}
\begin{abstract}
We construct 3D polytopes that admit a scaling with maximum likelihood (ML) degree one by using codimension zero toric 
fiber products of lower-dimensional polytopes and classify all the 3D polytopes that can be constructed in this way.
We also study a family of 3D polytopes whose faces admit ML degree one scalings.
We give a tight lower bound for the ML degree of these polytopes with a square base and we show that every value between this bound and
the normalized volume of the polytope can be attained.
These values depend on counting common and repeated roots associated to certain univariate
polynomials. This is similar in spirit to the description of the possible ML degrees of rational normal scrolls.
\end{abstract}

\maketitle

\section{Introduction}
Let $A\in \Z^{d\times n}$ be an integer matrix, the columns $a_1,\ldots,a_n$ of $A$ are lattice points in $\R^d$ and to each column we associate 
a scaling $c_i$. We define the \emph{scaled projective toric variety} $X_{A,c}$ to be the Zariski closure of the image of the
map $\Psi: (\C^{\ast})^d\to \pp^{n-1}$, $(t_1,\ldots,t_d)\mapsto (c_1t^{a_1}:\cdots:c_nt^{a_n})$, where $t^{a_i}:=t_1^{a_{1i}}t_2^{a_{2i}}\cdots t_d^{a_{di}}$.
We restrict to the case that the $c_i$ are positive real numbers and consider the map $\phi: X_{A,c}\setminus \HH \to \C^n$, $[p_1:\cdots:p_n]\mapsto \frac{1}{p_1+\cdots+p_n}(p_1,\ldots,p_n)$, where $\mathcal{H}=V(p_1p_2\cdots p_n(p_1+\ldots+p_n))$. The set $\MM_{A,c}=\Delta_{n-1}^{\circ}\cap \phi(X_{A,c}\setminus \HH)$
is the statistical model associated to $X_{A,c}$. In this setting, $\Delta_{n-1}^{\circ}$ denotes the open probability simplex, $\Delta_{n-1}^{\circ}=\{(p_1,\ldots,p_n):p_i\in \R, p_i>0, p_1+\cdots+p_n=1\}$. Points in $\Delta_{n-1}^{\circ}$ correspond to probability distributions of a discrete random variable with $n$ outcomes and a discrete
statistical model is simply a subset of $\Delta_{n-1}^{\circ}$ for some $n$. We are interested in the maximum likelihood estimation problem
for the class of models $\MM_{A,c}$, these are known as log-linear models in the statistical literature and also correspond to discrete (categorical) exponential
families. This class of models  is useful for the analysis of association in multivariate categorical data such as contingency tables \cite{bishop2007}.

The classical method to compute the MLE for a log-linear model $\MM_{A,c}$ is by using an Iterative Proportional Fitting (IPF) algorithm \cite{bishop2007}. The success of this algorithm to compute MLEs  relies on Birch's theorem \cite[Corollary 7.3.9] {sullivant2019} which states that for a positive data vector the MLE always exists and is unique.  The IPF algorithm has
several good properties such as guaranteed convergence to the unique MLE and convergence in one cycle of iteration whenever the MLE can
be written as a closed form expression of the data. 
There are results that identify classes of matrices $A$
for which IPF converges in one cycle \cite{coons24},  but no complete characterization of exactly when this happens.
An alternative algebraic geometric approach
to study the maximum likelihood estimation problem for the class of log-linear models is based on the \emph{maximum likelihood degree} of the model, also known as the ML degree.
The ML degree of an algebraic variety in $\pp^{n-1}$ is the number of complex critical points of the likelihood function, for generic data, that lie outside a specific hyperplane arrangement \cite{huhsturmfels2014}. The ML degree is an invariant of the embedded variety and it quantifies the algebraic complexity of obtaining the MLE
for the associated statistical model. If the variety under consideration has ML degree one, this implies that the maximum likelihood estimator is given by a closed form formula which is a rational function of the data \cite{huh2014,duarte2021}. 

In recent years, understanding the ML degrees of scaled projective toric varieties has been a central topic of interest in algebraic statistics. We say that the polytope $P=\conv(A)$ has ML degree $k$ if it admits a scaling such that $\mldeg(X_{A,c})=k$ and we adopt the convention that $A=P\cap \Z^d$, where the lattice points of $P$
affinely generate $\Z^d$ (see Section~\ref{sec:prelim}).
There is an important body of work that computes the ML degrees of different algebraic varieties arising from models in concrete 
applications such as hierarchical models and phylogenetics \cite{amendola2019,brysiewicz2023,garcia2023}.
Another research direction is to identify classes of scaled projective varieties whose ML degree is one \cite{coons2021,coons2024,davies2023,clarke2020,geiger2006}.
Finally a third and more recent trend focuses on understanding the possible ML degrees that arise by fixing $A$ and changing the choice $c$ of scalings \cite{clarke2024,AKK20,amendola2024}.
This leads to the study of the so-called Euler stratifications introduced by Telen and Wiesmann \cite{telen2024} and motivated by applications in physics.

\begin{figure}
    \centering
    \begin{tikzpicture}[line width =1.5pt,scale=0.5]
    \tikzstyle help lines=[color=blue!50,very thin]
    \draw (0,0) -- (5,0) -- (0,5) -- (0,0);
    \draw (10,0) -- (16,0) -- (16,5) -- (10,5) -- (10,0);
    \draw (21,0) -- (29,0) -- (23,5) -- (21,5) -- (21,0);
    \filldraw [black] (0,0) circle (4pt);
    \filldraw [black] (5,0) circle (4pt);
    \filldraw [black] (0,5) circle (4pt);
    \filldraw [black] (10,0) circle (4pt);
    \filldraw [black] (16,0) circle (4pt);
    \filldraw [black] (16,5) circle (4pt);
    \filldraw [black] (10,5) circle (4pt);
    \filldraw [black] (21,0) circle (4pt);
    \filldraw [black] (29,0) circle (4pt);
    \filldraw [black] (23,5) circle (4pt);
    \filldraw [black] (21,5) circle (4pt);

    \node at (2.5,-1) {$b\Delta_2$} ;
    \node at (13,-1) {$a\Delta_1\times b\Delta_1$} ;
    \node at (25,-1) {$T_{a,b,d}$} ;
\end{tikzpicture}

 \caption{Three families of polygons that admit scalings with ML degree one \cite{davies2023}.
 The vertex description of each family is $b\Delta_2=\conv((0,0),(0,b),(b,0))$, $a\Delta_1\times b\Delta_1=\conv((0,0),(a,0),(a,b),(0,b))$,  and $T_{a,b,d}=\conv((0,0),(a+db,0),(a,b),(0,b))$, where $a,d$ are nonnegative integers and $b$ is a positive integer.}
    \label{fig:2Dmldegree-one}
\end{figure}
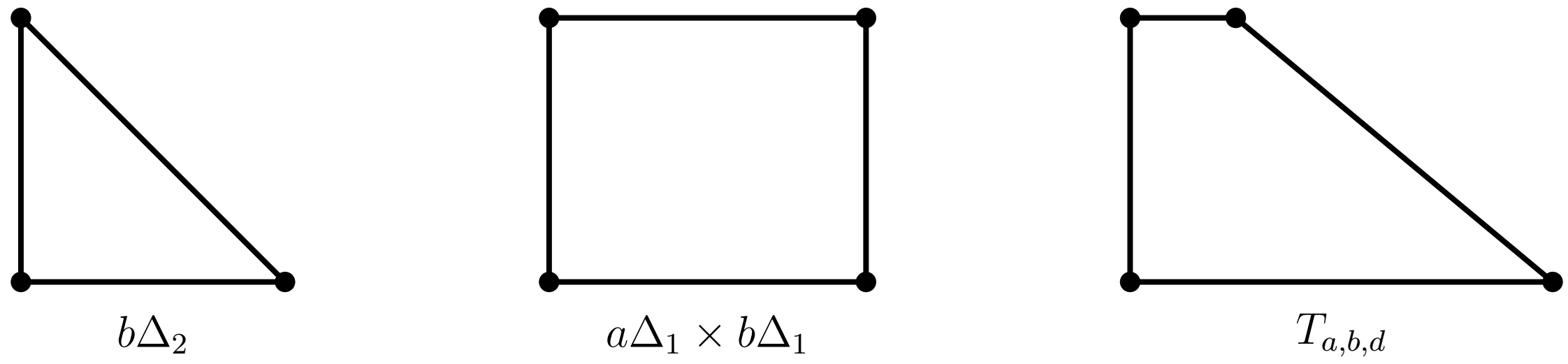

We are interested in studying the ML degree of sets of polytopes $P=\conv(A)$ in 3D that are built from 2D
models that have ML degree one. The motivation to study these polytopes comes from the following two results: first, the classification of 2D polytopes $P=\conv(A)$ with ML degree one \cite{davies2023,vonBothmer2010,garciapuente2010}. There are three infinite families of such polytopes; these are presented in Figure~\ref{fig:2Dmldegree-one}. In 3D and higher, the classification of polytopes for which there is a scaling such that the associated scaled projective toric variety has ML degree one is an open problem. 
The second is the result by Coons and Sullivant \cite{coons2021} that states that if the ML degree of $X_{A,c}$ is one and $P=\conv(A)$, then for any face $F$ of $P$ and $A' = F \cap A$, then $X_{A',c'}$ also has ML degree one, where $c'$ is the subvector of scalings restricted to the columns in $A'$. This result was generalized to higher ML degrees and states that $\mldeg(X_{A',c'})\leq \mldeg(X_{A,c})$ \cite{amendola2025}. Thus, to obtain a classification of 3D polytopes with ML degree one, it is reasonable to look into the properties and construction of 3D polytopes whose faces have ML degree one.
In addition to the ideas presented in the previous paragraph, it is important to note
that there are several 
known classes of polytopes, and associated statistical models, with ML degree one. These include decomposable models and certain classes of quasi-independence models. Even though these examples increase the understanding of which polytopes have ML degree one, it is difficult to gain insight from these classes because their dimension is high. To illustrate, consider the class of decomposable models. The first non trivial example with lowest dimension is the independence model of two binary random variables, this has dimension two and
corresponds
to the unit square. The next case to consider would be a decomposable model of three random variables. The full independence model corresponds to the unit cube and the next nontrivial example already has dimension four. Therefore, another reason to consider the construction
of 3D polytopes as we propose, is the fact that there are not many known examples of these beyond the classes considered in \cite[Section 3]{davies2023}.

In the first part of this article we construct 3D polytopes that admit a scaling with ML degree one by using codimension zero toric fiber products of 
lower-dimensional polytopes and classify all the 3D polytopes that can be constructed in this way.
In the second part, we study one family of 3D polytopes that we call tubular and are roughly obtained as the Cayley sum of two trapezoids, whose faces admit ML degree one scalings.
We give a tight lower bound for the ML degrees of polytopes in this family that have a square base (Definition~\ref{def:squarebase})  and we  prove that every ML degree between the lower bound and the normalized volume of the polytope can be attained.
We show that these values depend on counting common and repeated roots associated to certain univariate polynomials. This result is similar to the theorem that describes the possible ML degrees that arise for rational normal scrolls \cite{amendola2019} and Hirzebruch surfaces \cite{telen2024}.

The families of polytopes considered in this article intersect nontrivially with the
family of polytopes with ML degree one in \cite[Section~3]{davies2023}. The goal in
\cite{davies2023} is to show that the polytopes in that family admit a scaling with
ML degree one and to explain the connection to the geometry of the primitive collections of the normal fan of the polytope. Our approach in the first part is different: we show that the polytopes
we construct are codimension zero toric fiber products of lower-dimensional polytopes,
and we classify all the 3D polytopes that arise in this way. This gives a structural
explanation for why they have ML degree one, also for the polytopes that appear in both
families. In the second part, we go beyond ML degree one and study all the ML degrees
that tubular polytopes with a square base attain, including those that do not admit a
scaling with ML degree one.
\section{Preliminaries}\label{sec:prelim}
We start by defining the ML degree of a projective algebraic variety $V\subset \pp^{n-1}$. We use homogeneous coordinates $[p_1:\ldots:p_n]$ for points in $\pp^{n-1}$.
To obtain a statistical model from $V$, we consider the map $\phi:V\setminus \mathcal{H}\to \R^n$, $[p_1:\ldots:p_n]\mapsto \frac{1}{p_1+\cdots+p_n}(p_1,\ldots,p_n)$ where $\mathcal{H}=V(p_1\cdots p_n(p_1+\cdots+p_n))$. The set $\mathcal{M}=\Delta_{n-1}^{\circ}\cap \phi(V\setminus \HH)$ is a subset of the 
$(n-1)$-dimensional probability simplex. Each point $(p_1,\ldots,p_n)\in \MM$ is a probability distribution for a discrete random variable $Z$ with $n$ outcomes, namely $p_i=\mathrm{Prob}(Z=i)$ for $i=1,\ldots,n$.

To do parameter estimation  we take an independent identically distributed sample $\mathcal{D}$ from the random variable $Z$. This sample determines a data vector $u=(u_1,\ldots,u_n)$, where $u_i$ is the number of times the outcome $i$ appeared in $\mathcal{D}$. The \emph{likelihood function} is the rational function
\[
\LL(p|u)=\frac{p_1^{u_1}\cdots p_n^{u_n}}{(p_1+\cdots+p_n)^{u_1+\cdots+u_n}}.
\]
A \emph{maximum likelihood estimate} for the model $\MM$ given $u$ is a maximizer $\hat{p}$ of the likelihood function over $\mathcal{M}$. One way to find such maximizers is by computing the critical points of $\LL(p|u)$ over $V$. This leads to the definition of ML degree. Denote by $V_{\mathrm{reg}}$ the set of points in $V$ that are non-singular.

\begin{definition}
    The \emph{maximum likelihood degree} of $V$, denoted $\mldeg(V)$, is the
    number of complex critical points of the likelihood function
    $\LL(\cdot|u)$ on $V_{\mathrm{reg}}\setminus \HH$ for generic vectors $u$. We say that $V$ has rational MLE if $\mldeg(V)=1$.
\end{definition}

\begin{example}
This example is contained in \cite{huhsturmfels2014}.
Consider the quadric $V \subseteq \mathbb{P}^2$ associated with the polynomial
    \[
        f(p_0,p_1,p_2) = \det\begin{pmatrix}
            2p_0 & p_1\\
            p_1 & 2p_2
        \end{pmatrix} = 4p_0p_2-p_1^2.
    \]
    This is known as the \textit{Hardy-Weinberg curve}, and it can be alternatively defined parametrically by the map
    \begin{align*}
    s & \mapsto [p_0(s):p_1(s):p_2(s)] = [s^2: 2s(1-s): (1-s)^2].    
    \end{align*}
    Given data $u=(u_0,u_1,u_2)$, its associated likelihood function is
    \[
    \LL((p_0(s),p_1(s),p_2(s))|{u_0,u_1,u_2}) =  \frac{p_0(s)^{u_0} p_1(s)^{u_1} p_2(s)^{u_2}}{(p_0(s)+p_1(s)+p_2(s))^{u_0+u_1+u_2}} =
2^{u_1} s^{2u_0+u_1} (1-s)^{u_1+2u_2}.
    \]
    We can compute its ML degree by observing that the likelihood function has three critical points at $s=0,\, s=1$ and $s={(2u_0+u_1})/({2u_0+2u_1+2u_2})$. Since $p(0)$ and $p(1)$ are points of the hyperplane arrangement $\mathcal{H}=V(p_0p_1p_2(p_0+p_1+p_2))$, the Hardy-Weinberg curve has ML degree $1$.
    
    An alternative approach to calculating its ML degree involves the \textit{likelihood correspondence} of $V$. This is defined as the closure of $\{(p,u) \mid p \in V_{\mathrm{reg}} \setminus \mathcal{H},\,  p \text{ is critical for }\mathcal{L}(\cdot\mid u) \}$ in $\mathbb{P}^2 \times \mathbb{P}^2$. Using Proposition 3.5 of \cite{huhsturmfels2014}, the \textit{likelihood correspondence} of $V$ is given by the equations
    \begin{equation*}
\label{eq:detdet0}
{\det}
\begin{pmatrix}
2 p_0 & p_1 \\
p_1 & 2 p_2 
\end{pmatrix} \,\, = \,\,
{\det}
\begin{pmatrix}
p_1 + 2 p_2 &  2 p_0 + p_1 \\
 u_1 + 2 u_2 
 & 2 u_0 + u_1 
\end{pmatrix} 
\,\, = \,\, 0.
\end{equation*}
Saturating the ideal generated by them with respect to $p_0 + p_1 + p_2$, we obtain two other equations
\[
2 (u_1+2u_2)p_0 =  (2u_0+u_1)p_1, \qquad
  (u_1+2u_2)p_1 =	2 (2u_0+u_1)p_2 
\]
that determine a unique solution for $(p_0,p_1,p_2)$ for fixed data $u=(u_0,u_1,u_2)$, confirming that $\mldeg V = 1$.  
\end{example}

\begin{example}
Consider the \textit{binary four-cycle} graphical model $\mathcal{M}$. It is defined as the model obtained from the projective variety $V \subseteq \mathbb{P}^{15}$ given by the parametrization
\[
p_{ijk\ell} = a_{ij} b_{jk}c_{k\ell}d_{i\ell}, \qquad i, j, k, \ell \in \{0,1\}.
\]
This model does not have a rational MLE. As proven in \cite[Theorem~4.4]{geiger2006}, in the case of undirected graphical models, having rational MLE is equivalent to being a \textit{decomposable} model, which $\mathcal{M}$ is not. Instead, it is possible to calculate that the ML degree of $V$ is $13$ \cite[p.~1484]{geiger2006}.
\end{example}

The ML degree was introduced in \cite{catanese2006,hosten2005}, it is an invariant of an algebraic statistical model that measures the algebraic complexity of maximum likelihood estimation. 
We are interested in the ML degrees of scaled projective toric varieties, in statistics these correspond to discrete regular exponential families, as explained in \cite[Section 6.2]{sullivant2019}. 
\subsection{Maximum likelihood degree of toric varieties}
To study the ML degree of scaled projective toric varieties we use information coming from the torus action and the discrete geometry of $P=\conv(A)$.
Huh and Sturmfels \cite{huhsturmfels2014} express the ML degree of  $X_{A,c}$ in terms of the complement of a hypersurface in the torus $(\C^{\ast})^d$.
The pair $(A,c)$ defines the Laurent polynomial
\[
f(t)=c_1t^{a_1}+c_2t^{a_2}+\cdots+c_nt^{a_n}.
\] Throughout, we assume that the lattice points $A=P\cap\Z^d$ affinely generate
the lattice $\Z^d$, that is, the differences $a_i-a_j$ span $\Z^d$ as a group.
This guarantees that the map $(\C^\ast)^d\to X_{A,c}$, $t\mapsto(c_1t^{a_1}:\cdots:c_nt^{a_n})$,
is injective, so that $X_{A,c}\setminus\HH$ is isomorphic to the hypersurface complement
in the next theorem and $\deg X_{A,c}=\vol(P)$.
\begin{theorem}\cite[Theorem 3.2]{huhsturmfels2014} \label{thm:HuhSturmfels}
    The ML degree of the $d$-dimensional toric variety $X_{A,c}\subset \pp^{n-1}$ is equal to $(-1)^d$ times the Euler characteristic of the very affine variety
    \[
    X_{A,c}\setminus \mathcal{H}\simeq \{t\in (\C^{\ast})^d: f(t)\neq 0\}.
    \]
    For generic $c$, the ML degree agrees with the degree of $X_{A,c}$, which is the normalized volume $\vol(P)$ of the $d$-dimensional lattice
    polytope $P$.
\end{theorem}
Thus, for generic $c$, we expect $
\mldeg(X_{A,c})=\vol(P)$.
The next theorem refines the previous result and states that $\vol(P)$ is always an upper bound for $\mldeg(X_{A,c})$.
\begin{theorem}\cite[Theorem 2, Corollary 8]{amendola2019}
    Let $X_{A,c}$ be the scaled projective variety defined by the integer matrix $A\in \Z^{d\times n} $ and the vector $c\in(\C^{\ast})^n$.
    Then $\mldeg(X_{A,c})<deg(X_A)=\vol(P)$
    if and only if $c$ is in the vanishing locus of the principal $A$-determinant of the toric variety $X_A$.
\end{theorem}

\subsection{Toric fiber products}
Seth Sullivant defined the \textit{toric fiber product} in \cite{Sullivant2007}, we use his definition here. Given a positive integer $r \in \Z_{>0}$ and $s,\,t \in \Z^r_{>0}$ two vectors of positive integers, let
\[
\mathbb{K}[x] = \mathbb{K}[x^i_j \mid i \in [r],\, j \in [s_i]]
\quad
\text{and}
\quad
\mathbb{K}[y] = \mathbb{K}[y^i_j \mid i \in [r],\, j \in [t_i]]
\]
be two multigraded polynomial rings such that
\[
\deg(x^i_j) = \deg(y^i_k) = \mathbf{a}^i 
\]
for some set $\mathcal{A} = \{\mathbf{a}^1,\dots,\mathbf{a}^r\} \subseteq \Z^{\ell}$ of vectors. From now on, we also assume that there exists a vector $\omega \in \mathbb{Q}^{\ell}$ such that $\omega^T \mathbf{a}^i = 1$ for all $i$. This implies that homogeneous ideals in $\mathbb{K}[x]$ and $\mathbb{K}[y]$ are also homogeneous with respect to the standard grading.

\begin{definition}
    Let $I \subseteq \mathbb{K}[x]$ and $J \subseteq \mathbb{K}[y]$ be homogeneous ideals with respect to the multigrading induced by $\mathcal{A}$. The \textit{toric fiber product} (TFP) of $I$ and $J$, denoted by $I \times_{\mathcal{A}} J$, is the kernel of the ring homomorphism
    \begin{align*}
        \phi_{I,J}: \mathbb{K}[z^i_{jk}\mid i \in [r],\, j \in [s_i],\, k \in [t_i]] &\longrightarrow \left( \mathbb{K}[x]/I \right) \otimes_{\mathbb{K}} \left( \mathbb{K}[y]/J \right) \\
        z^i_{jk} &\longmapsto x^i_j \otimes y^i_k.
    \end{align*}
\end{definition}
An explicit description of the toric fiber product is given in \cite{engstrom2014} for the case of toric ideals.
For any configuration of integer vectors $\mathcal{P} \in \Z^{d \times r}$ , we denote by
\[
I_{\mathcal{P}} = (x^{\mathbf{u}}-x^{\mathbf{v}}\mid \mathcal{P}(u-v)=0)
\]
its corresponding toric ideal. Note that two configurations of integer vectors define the same ideal whenever they have the same kernel. Saying that $I_{\mathcal{P}}$ is homogeneous with respect to the multigrading induced by $\mathcal{A}$ is equivalent to requiring that $\mathcal{P} = \{\mathbf{p}^i_j \mid i \in [r],\,j \in [s_i]\}$ is partitioned so that there exists a linear map $\pi: \Z^d \to \Z^{\ell}$ which satisfies $\pi(\mathbf{p}^i_j)=\mathbf{a}^i$ for all $i$ and $j$. For brevity, by slight abuse of notation we may simply say that $\mathcal{P}$ is homogeneous with respect to $\mathcal{A}$.

\begin{proposition}
    Let $\mathcal{B} = \{\mathbf{b}^i_j \mid i \in [r],\,j \in [s_i]\}$ and $\mathcal{C} = \{\mathbf{c}^i_j \mid i \in [r],\,j \in [t_i]\}$ be two configurations of vectors that are homogeneous with respect to $\mathcal{A}$. Then, their toric fiber product is the toric ideal
    \[
    I_{\mathcal{B}} \times_{\mathcal{A}} I_{\mathcal{C}} = I_{\mathcal{B \times_{\mathcal{A}}\mathcal{C}}}
    \]
    where
    \[
    \mathcal{B \times_{\mathcal{A}}\mathcal{C}} = \{ (\mathbf{b}^i_j,\,\mathbf{c}^i_k)\mid i \in [r],\, j \in [s_i],\, k \in [t_i]\}.
    \]
\end{proposition}

It is of particular interest to study \textit{codimension zero toric fiber products}, which are TFPs in which $\mathcal{A}$ consists of linearly independent vectors. In this setting, it is possible to obtain information about the resulting ideal $I_{\mathcal{B \times_{\mathcal{A}}\mathcal{C}}}$ 
by using information from each of the factors $I_{\mathcal{B}},I_{\mathcal{C}}$. A relevant example is the multiplicative behaviour of the ML degree under the toric fiber product.

\begin{proposition}[\cite{AKK20}] \label{prop:mldeg_for_tfp}
    Let $\mathcal{B}$ and $\mathcal{C}$ be configurations of vectors that are homogeneous with respect to $\mathcal{A}$, and let $c_{\mathcal{B}}$ and $c_{\mathcal{C}}$ be vectors of weights for $\mathcal{B}$ and $\mathcal{C}$ respectively. If $\mathcal{A}$ consists of linearly independent vectors, then 
    \[
    \mldeg (X_{\mathcal{B}\times_\mathcal{A}\mathcal{C},\,(c_{\mathcal{B}}\times_{\mathcal{A}} c_{\mathcal{C}})}) = \mldeg (X_{\mathcal{B},\,c_{\mathcal{B}}}) \cdot \mldeg(X_{\mathcal{C},\,c_{\mathcal{C}}}),
    \]
    where $c_{\mathcal{B}}\times_{\mathcal{A}} c_{\mathcal{C}}$ is the vector of weight for $\mathcal{B}\times_\mathcal{A}\mathcal{C}$ such that $(c_{\mathcal{B}}\times_{\mathcal{A}} c_{\mathcal{C}})^i_{jk}:=(c_{\mathcal{B}})^i_j (c_{\mathcal{C}})^i_k$.
\end{proposition}
\begin{remark}
    The original statement of Proposition~\ref{prop:mldeg_for_tfp} in \cite{AKK20} uses weights that are all equal to one. However, the same proof applies using arbitrary weights for the vector configurations $\mathcal{B}$ and $\mathcal{C}$. 
\end{remark}

\begin{example}
    Let $r=3$ and $\mathcal{P} := \{p^1_1=(0\;0)^T,\,p^2_1=(1\;0)^T,\,p^2_2=(0\;1)^T,\,p^3_1=(1\;1)^T\}$ and $\pi(x,y) = x+y$, so $\mathcal{A} = \{
        0, 1, 2\} \subseteq \R$. Then $\mathcal{P} \times_{\mathcal{A}} \mathcal{P}$, up to an appropriate projection in the hyperplane $\{x_1 + x_2 = x_3 + x_4\}$ in which the TFP is contained, is the polytope
    
\[ \operatorname{conv}
    \begin{pmatrix}
        0 & 0 & 1 & 1 & 0 & 1 \\
        0 & 1 & 0 & 1 & 0 & 1 \\
        0 & 0 & 0 & 0 & 1 & -1 
    \end{pmatrix}.
    \]
    One can calculate that the variety $X_{\mathcal{P} \times_{\mathcal{A}} \mathcal{P}}$ has ML degree $4$. Since $X_\mathcal{P}$ has ML degree $1$, this example shows that Proposition \ref{prop:mldeg_for_tfp} does not hold if we remove the condition of $\mathcal{A}$ being a set of linearly independent vectors.
\begin{figure}[ht!]
    \centering
    \begin{tikzpicture}
    \begin{scope}[line width =1.5pt,scale=0.5]
    \tikzstyle help lines=[color=blue!50,very thin]
    \draw (0,0) -- (5,0) -- (5,5) -- (0,5) -- (0,0);
    \draw (10,0) -- (15,0) -- (15,5) -- (10,5) -- (10,0);
    \node at (7.5,2.5) {$\bigtimes_{\mathcal{A}}$};
    \node at (17.5, 2.5) {$\longrightarrow$};
    \filldraw [green] (0,0) circle (4pt);
    \filldraw [blue] (5,0) circle (4pt);
    \filldraw [blue] (0,5) circle (4pt);
    \filldraw [orange] (5,5) circle (4pt);
    \filldraw [green] (10,0) circle (4pt);
    \filldraw [blue] (15,0) circle (4pt);
    \filldraw [blue] (10,5) circle (4pt);
    \filldraw [orange] (15,5) circle (4pt);
    \end{scope}
    \tdplotsetmaincoords{70}{110} 
    \begin{scope}[tdplot_main_coords, xshift=11cm, scale=0.45, yshift=4.45cm]
            
        \coordinate (A) at (0,0,0);
        \filldraw [blue] (A) circle (4pt);        
        \coordinate (B) at (5,0,0);
        \filldraw [blue] (B) circle (4pt);    
        \coordinate (C) at (5,5,0);
        \filldraw [blue] (C) circle (4pt);    
        \coordinate (D) at (0,5,0);
        \filldraw [blue] (D) circle (4pt);    
        \coordinate (top) at (0,0,5);
        \filldraw [orange] (top) circle (4pt);    
        \coordinate (bottom) at (5,5,-5);
        \filldraw [green] (bottom) circle (4pt);    

        \draw[thick, dashed] (A) -- (B);
        \draw[thick] (B) -- (C) -- (D);
        \draw[thick, dashed] (D) -- (A);
        
        \draw[thick,dashed] (top) -- (A);
        \draw[thick] (top) -- (B);
        \draw[thick] (top) -- (C);
        \draw[thick] (top) -- (D);        

        \draw[thick,dashed] (bottom) -- (A);
        \draw[thick] (bottom) -- (B);
        \draw[thick] (bottom) -- (C);
        \draw[thick] (bottom) -- (D);   
        
    \end{scope}

\end{tikzpicture}
    \caption{Example of a TFP with linearly dependent $\mathcal{A}$ for which the multiplicative property of ML degrees fails.}
\end{figure}
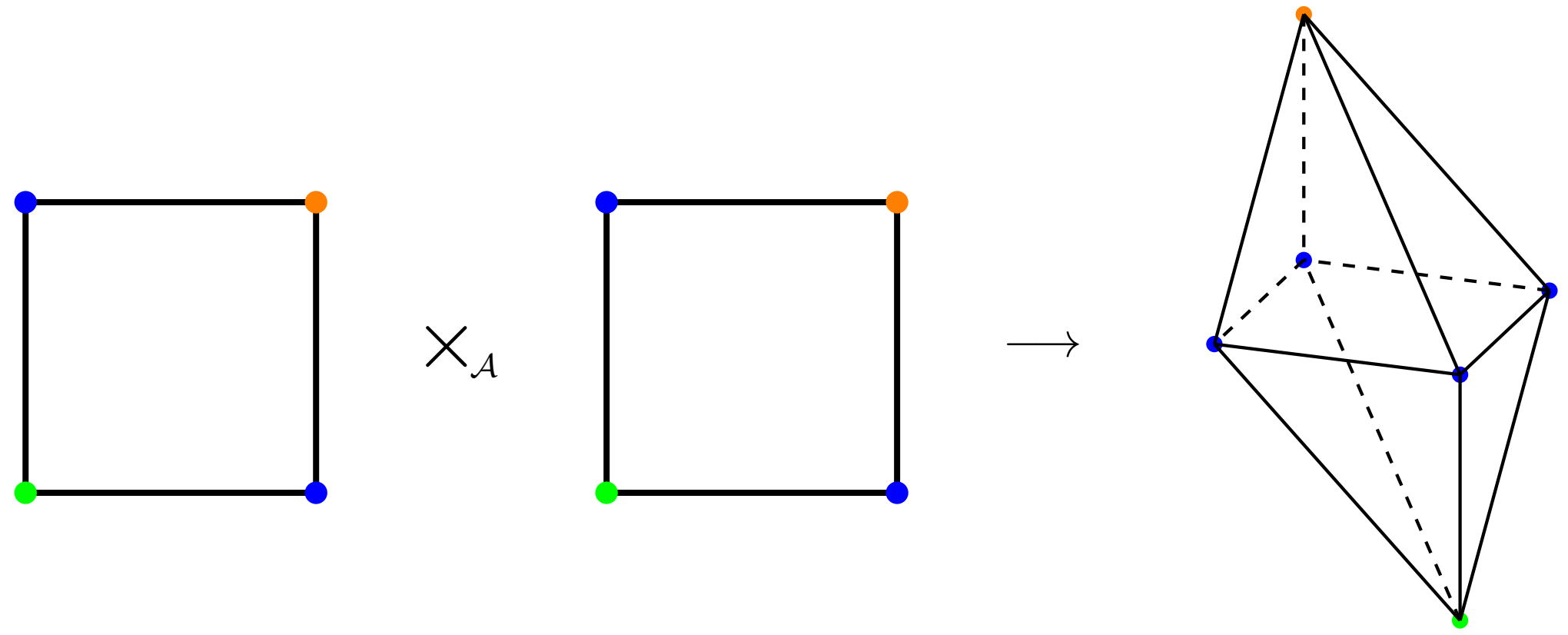
    
\end{example}

Given a lattice polytope $P$ with integer points $A \subseteq \Z^d$, we denote $I_P$ as the toric ideal associated to the matrix
\[
\hat{A} = \begin{pmatrix}
    1 & \cdots & 1 \\
     & A
\end{pmatrix}.
\]

\begin{remark}
    Saying that $I_{{P}}$ is homogeneous with respect to the multigrading induced by $\mathcal{A}$ is equivalent to requiring that $\mathcal{P} = \{\mathbf{p}^i_j \mid i \in [r],\,j \in [s_i]\}$ is partitioned so that there exists an {affine} map $\pi: \Z^d \to \Z^{\ell}$ which satisfies $\pi(\mathbf{p}^i_j)=\mathbf{a}^i$ for all $i$ and $j$. In fact, the restriction on the hyperplane $\{(x_0,\dots,x_n) \in \Z^{d+1} \mid x_0 = 1 \}$ gives a bijection between linear maps from $\Z^{d+1}$ to $\Z^{\ell}$ and affine maps from $\Z^d$ to $\Z^{\ell}$.
\end{remark}

\section{Toric fiber products of 2D polytopes}

In this section, we classify the $3$-dimensional polytopes with ML degree one that arise as codimension zero TFPs of lower dimensional polytopes. To obtain this classification, we first establish a relation between toric fiber products and Cayley sums of vector configurations. We then prove the classification in Theorem~\ref{thm:classification}.  
Our result also provides a complete description of all $3$-dimensional codimension zero TFPs of lower-dimensional polytopes, by Proposition~\ref{prop:3D_TFP}. 

The table in Figure~\ref{tb:polyhedra} lists the families of $3$-dimensional polytopes obtainable in this way. For every polytope in this list, there exists a positive scaling for which the ML degree of the associated scaled toric variety is $1$. The last column of the table gives the maximum ML degree that the polytope can attain, namely its normalized volume.
\newpage

\par\nointerlineskip
\begingroup
\setlength{\figureheight}{\dimexpr\textheight-\intextsep\relax}
\ifdim\pagegoal<\maxdimen
  \setlength{\figureheight}{%
    \dimexpr\pagegoal-\pagetotal-\pagedepth-2\intextsep\relax}
\fi

\begin{table}[H]
  \centering
  \tdplotsetmaincoords{70}{110}
  \tikzset{every picture/.style={baseline=(base.base),very thick,x=\baselinedim,y=\baselinedim}}

  \setbox\figurecaptionbox=\vbox{%
    \hsize=\linewidth
    \caption{Families of $3$-dimensional polytopes obtained as TFPs of polygons with ML degree one.}
    \label{tb:polyhedra}
  }
  \addtolength{\figureheight}{-\ht\figurecaptionbox}
  \addtolength{\figureheight}{-\dp\figurecaptionbox}

  \resizebox*{!}{\figureheight}{%
\begin{tabular}{c|c|c}
  Factors & Toric Fiber Product & ML degree for a general scaling \\ \hline
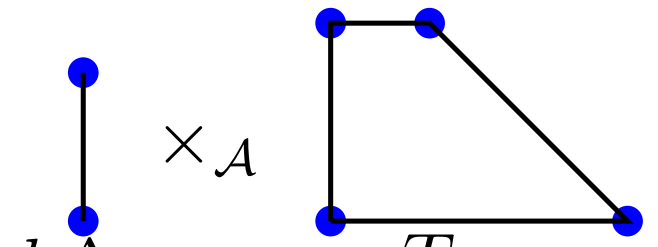
    \begin{tikzpicture}[line width =1.5pt,scale=1.12]
    \tikzstyle help lines=[color=blue!50,very thin]

    \node at (1.25,0.75) {$\times_{\mathcal{A}}$};
    \node at (0,-0.5) {$h\Delta_1$};
    \node at (4,-0.5) {$T_{a,b,d}$};
    \coordinate (A) at (0,0);
    \filldraw [blue] (A) circle (1.5pt);        
    \coordinate (B) at (0,1.5);
    \filldraw [blue] (B) circle (1.5pt);

    \draw[thick] (A) -- (B);
    
    \coordinate (A1) at (2.5,0);
    \filldraw [blue] (A1) circle (1.5pt);        
    \coordinate (B1) at (5.5,0);
    \filldraw [blue] (B1) circle (1.5pt);    
    \coordinate (C1) at (3.5,2);
    \filldraw [blue] (C1) circle (1.5pt);    
    \coordinate (D1) at (2.5,2);
    \filldraw [blue] (D1) circle (1.5pt); 

    \draw[thick] (A1) -- (B1) -- (C1) -- (D1) -- (A1);
    
    \addvmargin\addbase   
    \end{tikzpicture} & 
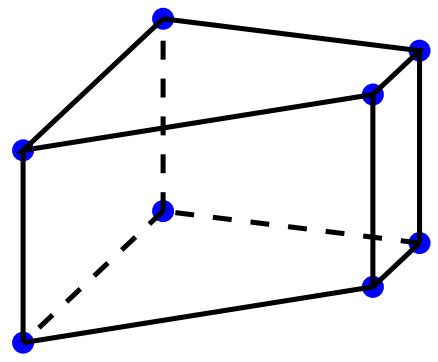
\begin{tikzpicture}[tdplot_main_coords,scale=0.76]
            
        \coordinate (A) at (0,0,0);
        \filldraw [blue] (A) circle (1.5pt);        
        \coordinate (B) at (3,0,0);
        \filldraw [blue] (B) circle (1.5pt);    
        \coordinate (C) at (1,2,0);
        \filldraw [blue] (C) circle (1.5pt);    
        \coordinate (D) at (0,2,0);
        \filldraw [blue] (D) circle (1.5pt);      
        \coordinate (A1) at (0,0,1.5);
        \filldraw [blue] (A1) circle (1.5pt);        
        \coordinate (B1) at (3,0,1.5);
        \filldraw [blue] (B1) circle (1.5pt);    
        \coordinate (C1) at (1,2,1.5);
        \filldraw [blue] (C1) circle (1.5pt);    
        \coordinate (D1) at (0,2,1.5);
        \filldraw [blue] (D1) circle (1.5pt);    

        \draw[thick, dashed] (A) -- (B);
        \draw[thick] (B) -- (C) -- (D);
        \draw[thick, dashed] (D) -- (A);

        \draw[thick] (A1) -- (B1) -- (C1) -- (D1) -- (A1);
        
        \draw[thick, dashed] (A) -- (A1);
        \draw[thick] (B) -- (B1);
        \draw[thick] (C) -- (C1);
        \draw[thick] (D) -- (D1);
        
      \addvmargin\addbase    
    \end{tikzpicture} & $h\left(6ab+3b^2d\right)$ \\ \hline
  \begin{tikzpicture}[line width =1.5pt,scale=1.12]
    \tikzstyle help lines=[color=blue!50,very thin]

    \node at (3.25,0.5) {$\times_{\mathcal{A}}$};
    \node at (1,-0.5) {$a\Delta_1 \times \Delta_1$};
    \node at (5,-0.5) {$a'\Delta_1 \times \Delta_1$};
    \coordinate (A) at (0,0);
    \filldraw [blue] (A) circle (1.5pt);        
    \coordinate (B) at (2,0);
    \filldraw [blue] (B) circle (1.5pt);    
    \coordinate (C) at (2,1);
    \filldraw [red] (C) circle (1.5pt);    
    \coordinate (D) at (0,1);
    \filldraw [red] (D) circle (1.5pt); 

    \draw[thick] (A) -- (B) -- (C) -- (D) -- (A);
    
    \coordinate (A1) at (4.5,0);
    \filldraw [blue] (A1) circle (1.5pt);        
    \coordinate (B1) at (5.5,0);
    \filldraw [blue] (B1) circle (1.5pt);    
    \coordinate (C1) at (5.5,1);
    \filldraw [red] (C1) circle (1.5pt);    
    \coordinate (D1) at (4.5,1);
    \filldraw [red] (D1) circle (1.5pt); 

    \draw[thick] (A1) -- (B1) -- (C1) -- (D1) -- (A1);
    
    \addvmargin\addbase   
    \end{tikzpicture} & 
    \begin{tikzpicture}[tdplot_main_coords,scale=0.76]
            
        \coordinate (A) at (0,0,0);
        \filldraw [blue] (A) circle (1.5pt);        
        \coordinate (B) at (1,0,0);
        \filldraw [blue] (B) circle (1.5pt);    
        \coordinate (C) at (1,2,0);
        \filldraw [blue] (C) circle (1.5pt);    
        \coordinate (D) at (0,2,0);
        \filldraw [blue] (D) circle (1.5pt);    
        \coordinate (top) at (0,0,1);   
        \coordinate (A1) at (0,0,1);
        \filldraw [red] (A1) circle (1.5pt);        
        \coordinate (B1) at (1,0,1);
        \filldraw [red] (B1) circle (1.5pt);    
        \coordinate (C1) at (1,2,1);
        \filldraw [red] (C1) circle (1.5pt);    
        \coordinate (D1) at (0,2,1);
        \filldraw [red] (D1) circle (1.5pt);    

        \draw[thick, dashed] (A) -- (B);
        \draw[thick] (B) -- (C) -- (D);
        \draw[thick, dashed] (D) -- (A);

        \draw[thick] (A1) -- (B1) -- (C1) -- (D1) -- (A1);
        
        \draw[thick, dashed] (A) -- (A1);
        \draw[thick] (B) -- (B1);
        \draw[thick] (C) -- (C1);
        \draw[thick] (D) -- (D1);
        
      \addvmargin\addbase    
    \end{tikzpicture}
  & $6aa'$ \\ \hline
  \begin{tikzpicture}[line width =1.5pt,scale=1.12]
    \tikzstyle help lines=[color=blue!50,very thin]

    \node at (3.25,0.5) {$\times_{\mathcal{A}}$};
    \node at (1,-0.5) {$a\Delta_1 \times \Delta_1$};
    \node at (5,-0.5) {$a'\Delta_1*\{0\}$};
    \coordinate (A) at (0,0);
    \filldraw [blue] (A) circle (1.5pt);        
    \coordinate (B) at (2,0);
    \filldraw [blue] (B) circle (1.5pt);    
    \coordinate (C) at (2,1);
    \filldraw [red] (C) circle (1.5pt);    
    \coordinate (D) at (0,1);
    \filldraw [red] (D) circle (1.5pt); 

    \draw[thick] (A) -- (B) -- (C) -- (D) -- (A);
    
    \coordinate (A1) at (4.5,0);
    \filldraw [blue] (A1) circle (1.5pt);        
    \coordinate (B1) at (5.5,0);
    \filldraw [blue] (B1) circle (1.5pt);      
    \coordinate (D1) at (4.5,1);
    \filldraw [red] (D1) circle (1.5pt); 

    \draw[thick] (A1) -- (B1) -- (D1) -- (A1);
    
    \addvmargin\addbase   
    \end{tikzpicture} & 
    \begin{tikzpicture}[tdplot_main_coords,scale=0.76]
            
        \coordinate (A) at (0,0,0);
        \filldraw [blue] (A) circle (1.5pt);        
        \coordinate (B) at (1,0,0);
        \filldraw [blue] (B) circle (1.5pt);    
        \coordinate (C) at (1,2,0);
        \filldraw [blue] (C) circle (1.5pt);    
        \coordinate (D) at (0,2,0);
        \filldraw [blue] (D) circle (1.5pt);    
        \coordinate (top) at (0,0,1);   
        \coordinate (A1) at (0,0,1);
        \filldraw [red] (A1) circle (1.5pt);        
        \coordinate (D1) at (0,2,1);
        \filldraw [red] (D1) circle (1.5pt);

        \draw[thick, dashed] (A) -- (B);
        \draw[thick] (B) -- (C) -- (D);
        \draw[thick, dashed] (D) -- (A);

        \draw[thick] (A1) -- (D1);
        
        \draw[thick, dashed] (A) -- (A1);
        \draw[thick] (D) -- (D1);
        \draw[thick] (B) -- (A1);
        \draw[thick] (C) -- (D1);
        
      \addvmargin\addbase    
    \end{tikzpicture}
  & $3aa'$ \\ \hline
  \begin{tikzpicture}[line width =1.5pt,scale=1.12]
    \tikzstyle help lines=[color=blue!50,very thin]

    \node at (3.25,0.5) {$\times_{\mathcal{A}}$};
    \node at (1,-0.5) {$a\Delta_1 \times \Delta_1$};
    \node at (5.5,-0.5) {$T_{a',1,d'}$};
    \coordinate (A) at (0,0);
    \filldraw [blue] (A) circle (1.5pt);        
    \coordinate (B) at (2,0);
    \filldraw [blue] (B) circle (1.5pt);    
    \coordinate (C) at (2,1);
    \filldraw [red] (C) circle (1.5pt);    
    \coordinate (D) at (0,1);
    \filldraw [red] (D) circle (1.5pt); 

    \draw[thick] (A) -- (B) -- (C) -- (D) -- (A);
    
    \coordinate (A1) at (4.5,0);
    \filldraw [blue] (A1) circle (1.5pt);        
    \coordinate (B1) at (6.5,0);
    \filldraw [blue] (B1) circle (1.5pt);    
    \coordinate (C1) at (5.5,1);
    \filldraw [red] (C1) circle (1.5pt);    
    \coordinate (D1) at (4.5,1);
    \filldraw [red] (D1) circle (1.5pt); 

    \draw[thick] (A1) -- (B1) -- (C1) -- (D1) -- (A1);
    
    \addvmargin\addbase   
    \end{tikzpicture} & \begin{tikzpicture}[tdplot_main_coords,scale=0.76]
            
        \coordinate (A) at (0,0,0);
        \filldraw [blue] (A) circle (1.5pt);        
        \coordinate (B) at (2,0,0);
        \filldraw [blue] (B) circle (1.5pt);    
        \coordinate (C) at (2,2,0);
        \filldraw [blue] (C) circle (1.5pt);    
        \coordinate (D) at (0,2,0);
        \filldraw [blue] (D) circle (1.5pt);    
        \coordinate (top) at (0,0,1);   
        \coordinate (A1) at (0,0,1);
        \filldraw [red] (A1) circle (1.5pt);        
        \coordinate (B1) at (1,0,1);
        \filldraw [red] (B1) circle (1.5pt);    
        \coordinate (C1) at (1,2,1);
        \filldraw [red] (C1) circle (1.5pt);    
        \coordinate (D1) at (0,2,1);
        \filldraw [red] (D1) circle (1.5pt);    

        \draw[thick, dashed] (A) -- (B);
        \draw[thick] (B) -- (C) -- (D);
        \draw[thick, dashed] (D) -- (A);

        \draw[thick] (A1) -- (B1) -- (C1) -- (D1) -- (A1);
        
        \draw[thick, dashed] (A) -- (A1);
        \draw[thick] (B) -- (B1);
        \draw[thick] (C) -- (C1);
        \draw[thick] (D) -- (D1);
        
      \addvmargin\addbase    
    \end{tikzpicture}
  & $a(6a'+3d')$ \\ \hline
\begin{tikzpicture}[line width =1.5pt,scale=1.12]
    \tikzstyle help lines=[color=blue!50,very thin]

    \node at (3.25,0.5) {$\times_{\mathcal{A}}$};
    \node at (1,-0.5) {$a\Delta_1*\{0\}$};
    \node at (5,-0.5) {$a'\Delta_1*\{0\}$};
    \coordinate (A) at (0,0);
    \filldraw [blue] (A) circle (1.5pt);        
    \coordinate (B) at (2,0);
    \filldraw [blue] (B) circle (1.5pt);       
    \coordinate (D) at (0,1);
    \filldraw [red] (D) circle (1.5pt); 

    \draw[thick] (A) -- (B) -- (D) -- (A);
    
    \coordinate (A1) at (4.5,0);
    \filldraw [blue] (A1) circle (1.5pt);        
    \coordinate (B1) at (5.5,0);
    \filldraw [blue] (B1) circle (1.5pt);        
    \coordinate (D1) at (4.5,1);
    \filldraw [red] (D1) circle (1.5pt); 

    \draw[thick] (A1) -- (B1) -- (D1) -- (A1);
    
    \addvmargin\addbase   
    \end{tikzpicture} & \begin{tikzpicture}[tdplot_main_coords,scale=0.76]
            
        \coordinate (A) at (0,0,0);
        \filldraw [blue] (A) circle (1.5pt);        
        \coordinate (B) at (2,0,0);
        \filldraw [blue] (B) circle (1.5pt);    
        \coordinate (C) at (2,2,0);
        \filldraw [blue] (C) circle (1.5pt);    
        \coordinate (D) at (0,2,0);
        \filldraw [blue] (D) circle (1.5pt);    
        \coordinate (top) at (0,0,1);   
        \filldraw [red] (top) circle (1.5pt);       

        \draw[thick, dashed] (A) -- (B);
        \draw[thick] (B) -- (C) -- (D);
        \draw[thick, dashed] (D) -- (A);
        
        \draw[thick, dashed] (top) -- (A);
        \draw[thick] (top) -- (B);
        \draw[thick] (top) -- (C);
        \draw[thick] (top) -- (D);
        
      \addvmargin\addbase    
    \end{tikzpicture}
  & $2aa'$ \\ \hline
  \begin{tikzpicture}[line width =1.5pt,scale=1.12]
    \tikzstyle help lines=[color=blue!50,very thin]

    \node at (3.25,0.5) {$\times_{\mathcal{A}}$};
    \node at (1,-0.5) {$a\Delta_1*\{0\}$};
    \node at (5,-0.5) {$a'\Delta_1*\{0\}$};
    \coordinate (A) at (0,0);
    \filldraw [red] (A) circle (1.5pt);        
    \coordinate (B) at (2,0);
    \filldraw [red] (B) circle (1.5pt);       
    \coordinate (D) at (0,1);
    \filldraw [blue] (D) circle (1.5pt); 

    \draw[thick] (A) -- (B) -- (D) -- (A);
    
    \coordinate (A1) at (4.5,0);
    \filldraw [blue] (A1) circle (1.5pt);        
    \coordinate (B1) at (5.5,0);
    \filldraw [blue] (B1) circle (1.5pt);        
    \coordinate (D1) at (4.5,1);
    \filldraw [red] (D1) circle (1.5pt); 

    \draw[thick] (A1) -- (B1) -- (D1) -- (A1);
    
    \addvmargin\addbase   
    \end{tikzpicture} & \begin{tikzpicture}[tdplot_main_coords,scale=0.76]
            
        \coordinate (A) at (0,0,0);
        \filldraw [blue] (A) circle (1.5pt);        
        \coordinate (B) at (1,0,0);
        \filldraw [blue] (B) circle (1.5pt);    
        \coordinate (A1) at (0,0,1);
        \filldraw [red] (A1) circle (1.5pt);        
        \coordinate (D1) at (0,2,1);
        \filldraw [red] (D1) circle (1.5pt);    

        \draw[thick, dashed] (A1) -- (A) -- (B);
        \draw[thick] (B) -- (A1) -- (D1) -- (B);
        \draw[thick, dashed] (D1) -- (A);
        
      \addvmargin\addbase    
    \end{tikzpicture}
  & $aa'$ \\ \hline
  \begin{tikzpicture}[line width =1.5pt,scale=1.12]
    \tikzstyle help lines=[color=blue!50,very thin]

    \node at (3.25,0.5) {$\times_{\mathcal{A}}$};
    \node at (1,-0.5) {$a\Delta_1 * \{0\}$};
    \node at (5.5,-0.5) {$T_{a',1,d'}$};
    \coordinate (A) at (0,0);
    \filldraw [blue] (A) circle (1.5pt);        
    \coordinate (B) at (2,0);
    \filldraw [blue] (B) circle (1.5pt);      
    \coordinate (D) at (0,1);
    \filldraw [red] (D) circle (1.5pt); 

    \draw[thick] (A) -- (B) -- (D) -- (A);
    
    \coordinate (A1) at (4.5,0);
    \filldraw [blue] (A1) circle (1.5pt);        
    \coordinate (B1) at (6.5,0);
    \filldraw [blue] (B1) circle (1.5pt);    
    \coordinate (C1) at (5.5,1);
    \filldraw [red] (C1) circle (1.5pt);    
    \coordinate (D1) at (4.5,1);
    \filldraw [red] (D1) circle (1.5pt); 

    \draw[thick] (A1) -- (B1) -- (C1) -- (D1) -- (A1);
    
    \addvmargin\addbase   
    \end{tikzpicture} & \begin{tikzpicture}[tdplot_main_coords,scale=0.76]
            
        \coordinate (A) at (0,0,0);
        \filldraw [blue] (A) circle (1.5pt);        
        \coordinate (B) at (2,0,0);
        \filldraw [blue] (B) circle (1.5pt);    
        \coordinate (C) at (2,2,0);
        \filldraw [blue] (C) circle (1.5pt);    
        \coordinate (D) at (0,2,0);
        \filldraw [blue] (D) circle (1.5pt);      
        \coordinate (A1) at (0,0,1);
        \filldraw [red] (A1) circle (1.5pt);        
        \coordinate (B1) at (1,0,1);
        \filldraw [red] (B1) circle (1.5pt);      

        \draw[thick, dashed] (A) -- (B);
        \draw[thick] (B) -- (C) -- (D);
        \draw[thick, dashed] (D) -- (A);

        \draw[thick] (A1) -- (B1);
        
        \draw[thick, dashed] (A) -- (A1);
        \draw[thick] (B) -- (B1);
        \draw[thick] (C) -- (B1);
        \draw[thick] (D) -- (A1);
        
      \addvmargin\addbase    
    \end{tikzpicture}
  & $a(3a'+2d')$ \\ \hline
    \begin{tikzpicture}[line width =1.5pt,scale=1.12]
    \tikzstyle help lines=[color=blue!50,very thin]

    \node at (3.25,0.5) {$\times_{\mathcal{A}}$};
    \node at (1,-0.5) {$a\Delta_1 \ast \{0\}$};
    \node at (5.5,-0.5) {$T_{a',1,d'}$};
    \coordinate (A) at (0,0);
    \filldraw [blue] (A) circle (1.5pt);        
    \coordinate (B) at (2,0);
    \filldraw [blue] (B) circle (1.5pt);       
    \coordinate (D) at (0,1);
    \filldraw [red] (D) circle (1.5pt); 

    \draw[thick] (A) -- (B) -- (D) -- (A);
    
    \coordinate (A1) at (4.5,0);
    \filldraw [red] (A1) circle (1.5pt);        
    \coordinate (B1) at (6.5,0);
    \filldraw [red] (B1) circle (1.5pt);    
    \coordinate (C1) at (5.5,1);
    \filldraw [blue] (C1) circle (1.5pt);    
    \coordinate (D1) at (4.5,1);
    \filldraw [blue] (D1) circle (1.5pt); 

    \draw[thick] (A1) -- (B1) -- (C1) -- (D1) -- (A1);
    
    \addvmargin\addbase   
    \end{tikzpicture} & \begin{tikzpicture}[tdplot_main_coords,scale=0.76]
            
        \coordinate (A) at (0,0,0);
        \filldraw [blue] (A) circle (1.5pt);        
        \coordinate (B) at (1,0,0);
        \filldraw [blue] (B) circle (1.5pt);    
        \coordinate (C) at (1,2,0);
        \filldraw [blue] (C) circle (1.5pt);    
        \coordinate (D) at (0,2,0);
        \filldraw [blue] (D) circle (1.5pt);      
        \coordinate (A1) at (0,0,1);
        \filldraw [red] (A1) circle (1.5pt);        
        \coordinate (B1) at (2,0,1);
        \filldraw [red] (B1) circle (1.5pt);      

        \draw[thick, dashed] (A) -- (B);
        \draw[thick] (B) -- (C) -- (D);
        \draw[thick, dashed] (D) -- (A);

        \draw[thick] (A1) -- (B1);
        
        \draw[thick, dashed] (A) -- (A1);
        \draw[thick] (B) -- (B1);
        \draw[thick] (C) -- (B1);
        \draw[thick] (D) -- (A1);
        
      \addvmargin\addbase    
    \end{tikzpicture}
  & $a(3a'+d')$ \\ \hline
  \begin{tikzpicture}[line width =1.5pt,scale=1.12]
    \tikzstyle help lines=[color=blue!50,very thin]

    \node at (3.75,0.5) {$\times_{\mathcal{A}}$};
    \node at (1.5,-0.5) {$T_{a,1,d}$};
    \node at (5.5,-0.5) {$T_{a',1,d'}$};
    \coordinate (A) at (0,0);
    \filldraw [blue] (A) circle (1.5pt);        
    \coordinate (B) at (3,0);
    \filldraw [blue] (B) circle (1.5pt);    
    \coordinate (C) at (2,1);
    \filldraw [red] (C) circle (1.5pt);    
    \coordinate (D) at (0,1);
    \filldraw [red] (D) circle (1.5pt); 

    \draw[thick] (A) -- (B) -- (C) -- (D) -- (A);
    
    \coordinate (A1) at (4.5,0);
    \filldraw [blue] (A1) circle (1.5pt);        
    \coordinate (B1) at (6.5,0);
    \filldraw [blue] (B1) circle (1.5pt);    
    \coordinate (C1) at (5.5,1);
    \filldraw [red] (C1) circle (1.5pt);    
    \coordinate (D1) at (4.5,1);
    \filldraw [red] (D1) circle (1.5pt); 

    \draw[thick] (A1) -- (B1) -- (C1) -- (D1) -- (A1);
    
    \addvmargin\addbase   
    \end{tikzpicture}  & \begin{tikzpicture}[tdplot_main_coords,scale=0.76]
            
        \coordinate (A) at (0,0,0);
        \filldraw [blue] (A) circle (1.5pt);        
        \coordinate (B) at (2,0,0);
        \filldraw [blue] (B) circle (1.5pt);    
        \coordinate (C) at (2,3,0);
        \filldraw [blue] (C) circle (1.5pt);    
        \coordinate (D) at (0,3,0);
        \filldraw [blue] (D) circle (1.5pt);    
        \coordinate (top) at (0,0,1);   
        \coordinate (A1) at (0,0,1);
        \filldraw [red] (A1) circle (1.5pt);        
        \coordinate (B1) at (1,0,1);
        \filldraw [red] (B1) circle (1.5pt);    
        \coordinate (C1) at (1,2,1);
        \filldraw [red] (C1) circle (1.5pt);    
        \coordinate (D1) at (0,2,1);
        \filldraw [red] (D1) circle (1.5pt);     

        \draw[thick, dashed] (A) -- (B);
        \draw[thick] (B) -- (C) -- (D);
        \draw[thick, dashed] (D) -- (A);

        \draw[thick] (A1) -- (B1) -- (C1) -- (D1) -- (A1);
        
        \draw[thick, dashed] (A) -- (A1);
        \draw[thick] (B) -- (B1);
        \draw[thick] (C) -- (C1);
        \draw[thick] (D) -- (D1);
        
      \addvmargin\addbase    
    \end{tikzpicture}
  & $6aa' + 3(ad' + a'd) + 2dd'$ \\ \hline
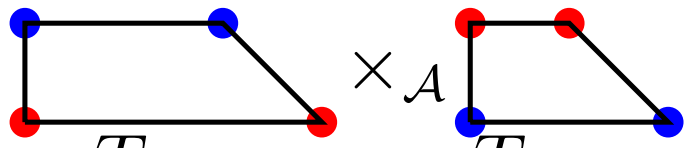
  \begin{tikzpicture}[line width =1.5pt,scale=1.12]
    \tikzstyle help lines=[color=blue!50,very thin]

    \node at (3.75,0.5) {$\times_{\mathcal{A}}$};
    \node at (1.5,-0.5) {$T_{a,1,d}$};
    \node at (5.5,-0.5) {$T_{a',1,d'}$};
    \coordinate (A) at (0,0);
    \filldraw [red] (A) circle (1.5pt);        
    \coordinate (B) at (3,0);
    \filldraw [red] (B) circle (1.5pt);    
    \coordinate (C) at (2,1);
    \filldraw [blue] (C) circle (1.5pt);    
    \coordinate (D) at (0,1);
    \filldraw [blue] (D) circle (1.5pt); 

    \draw[thick] (A) -- (B) -- (C) -- (D) -- (A);
    
    \coordinate (A1) at (4.5,0);
    \filldraw [blue] (A1) circle (1.5pt);        
    \coordinate (B1) at (6.5,0);
    \filldraw [blue] (B1) circle (1.5pt);    
    \coordinate (C1) at (5.5,1);
    \filldraw [red] (C1) circle (1.5pt);    
    \coordinate (D1) at (4.5,1);
    \filldraw [red] (D1) circle (1.5pt); 

    \draw[thick] (A1) -- (B1) -- (C1) -- (D1) -- (A1);
    
    \addvmargin\addbase   
    \end{tikzpicture}  & \begin{tikzpicture}[tdplot_main_coords,scale=0.76]
            
        \coordinate (A) at (0,0,0);
        \filldraw [blue] (A) circle (1.5pt);        
        \coordinate (B) at (2,0,0);
        \filldraw [blue] (B) circle (1.5pt);    
        \coordinate (C) at (2,1,0);
        \filldraw [blue] (C) circle (1.5pt);    
        \coordinate (D) at (0,1,0);
        \filldraw [blue] (D) circle (1.5pt);    
        \coordinate (top) at (0,0,1);   
        \coordinate (A1) at (0,0,1);
        \filldraw [red] (A1) circle (1.5pt);        
        \coordinate (B1) at (1,0,1);
        \filldraw [red] (B1) circle (1.5pt);    
        \coordinate (C1) at (1,3,1);
        \filldraw [red] (C1) circle (1.5pt);    
        \coordinate (D1) at (0,3,1);
        \filldraw [red] (D1) circle (1.5pt);    

        \draw[thick, dashed] (A) -- (B);
        \draw[thick] (B) -- (C);
        \draw[thick, dashed] (C) -- (D) -- (A);

        \draw[thick] (A1) -- (B1) -- (C1) -- (D1) -- (A1);
        
        \draw[thick, dashed] (A) -- (A1);
        \draw[thick] (B) -- (B1);
        \draw[thick] (C) -- (C1);
        \draw[thick, dashed] (D) -- (D1);
        
      \addvmargin\addbase    
    \end{tikzpicture}
  & $6aa' + 3(ad'+a'd) + dd'$ 

\end{tabular}
  }\par
  \nointerlineskip
  \usebox{\figurecaptionbox}
\end{table}
\endgroup

\subsection{Cayley sum and linear independence of $\mathcal{A}$}
Let $\mathcal{A} = \{\mathbf{a}^1,\dots,\mathbf{a}^r\} \subseteq \Z^{\ell}$ be a set of vectors. Throughout this section, we assume that $\mathcal{A}$ consists of linearly independent vectors. We also denote by $P$ and $Q$ two lattice polytopes, with respective sets of lattice points $A$ and $B$.

\begin{definition}
    Let $\mathcal{P}_0,\dots,\mathcal{P}_m \subseteq \Z^s$ be  finite configurations of vectors. The \textit{Cayley sum} $\mathcal{P}_0 * \mathcal{P}_1 *\dots * \mathcal{P}_m$ is defined as
    \[
    \left( \mathcal{P}_0 \times \{0\} \right) \cup \left( \mathcal{P}_1 \times \{e_1\} \right) \cup \dots \cup \left( \mathcal{P}_m \times \{e_m\} \right) \subseteq \Z^s \times \Z^m
    \]
    where $e_1,\dots,e_m$ is the standard basis of $\Z^m$.
\end{definition}
In the following discussion, whenever we have an affine map $f: \Z^{d_1}\to\Z^{d_2}$ we denote by $f_{\R}:\R^{d_1}\to\R^{d_2}$ the affine map obtained by tensoring $f$ with $\R$. 
\begin{lemma} \label{lemma:CayleySum}
  If $P \subseteq \R^d$ is full-dimensional and $A$ is homogeneous with respect to $\mathcal{A} = \{\mathbf{a}^1,\dots,\mathbf{a}^r\}$ via an affine map $\pi: \Z^d \to \Z^{\ell}$ of rank $k$, then 
    \begin{enumerate}
        \item we have $r=k+1$;
        \item the points of $A$ lie in $r$ parallel affine subspaces of codimension $k$;
        \item the polytopes $P_1,\dots,P_r$ defined as 
        \[
        P_i := P \cap \pi^{-1}_{\R}(\mathbf{a}^i) \quad\forall i=1,\dots,r ,
        \]
        are faces of $P$;
        \item the toric ideal associated to $P$ is
        \[
        I_{P} = I_{P_1*\dots*P_r}.
        \]
    \end{enumerate}
\end{lemma}
\begin{proof}
    \begin{enumerate} 
        \item Since $P$ is full-dimensional, the affine span of $\pi(A)  = \mathcal{A}$ has dimension $k$. Since $\mathcal{A}$ is made of affinely independent vectors, its image is  a simplex of dimension $r-1$, thus $r = k+1$. 
\item The set $A$ is contained in the union of the fibers $\cup_{i=1}^r\pi^{-1}_{\R}(\mathbf{a}^i)$. Each of the fibers is an affine subspace of codimension $k$ that is parallel to $\ker (\pi_{\R} - \pi_{\R}(0))$. 
        \item Since affine maps commute with the convex hull operator, the image of $P$ by $\pi_{\R}$ is the polytope $P':= \pi_{\mathbb{R}}(\conv(A)) =\conv(\pi_{\mathbb{R}}(A)) = \conv(\mathbf{a}^1,\dots,\mathbf{a}^r)$. By hypothesis, the vectors $\mathbf{a}^i$ are affinely independent, hence they are precisely the vertices of $P'$. Let $H_1,\dots,H_r \subseteq \R^\ell$ be supporting hyperplanes that meet the polytope at the vertices $\mathbf{a}^1,\dots,\mathbf{a}^r$ respectively, then $\pi^{-1}_{\R}(H_1),\dots,\pi^{-1}_{\R}(H_{r}) \subseteq \R^d$ are the supporting hyperplanes of $P_1,\dots,P_r$ respectively. 
        \item Since $\ker (\pi-\pi(0)) \subseteq \Z^d$ is a primitive sublattice, there exists a sublattice $W \subseteq \Z^d$ such that $\Z^d = \ker (\pi-\pi(0)) \oplus W$. Up to a unimodular transformation, we can assume that $\ker (\pi-\pi(0)) = \Z^{d-k} \times \{0\}^k \subseteq \Z^d$ and $W= \{0\}^{d-k} \times \Z^k \subseteq \Z^d$. Without loss of generality, we can also assume that the preimage of $\mathcal{A}$ under $\pi$ is of the form $\Z^{d-k}\times \{0,w_1,\dots,w_{r-1}\}$ for some linearly independent vectors $w_1,\dots,w_{r-1} \in W$. Indeed, translations do not change the toric ideal associated to a polytope. Notice also that $P_i \cap \Z^d \subseteq \Z^{d-k} \times \{w_{i-1}\}$ are lattice polytopes with integer points $A_i \times \{w_{i-1}\}$. By choosing an invertible matrix $U \in \mathbb{Q}^{k \times k}$ such that $Uw_i = e_i \subseteq \Z^k$, we obtain 
        \[
        (\mathrm{Id}_{{d-k}} \oplus U) A = A_1*\dots*A_r,
        \]
        which implies that $A$ and $A_1*\dots*A_r$ have the same kernel. \qedhere
    \end{enumerate}
\end{proof}

\begin{remark} \label{remark:simplex}
    By the previous lemma, if $A \subseteq \Z^d$ is homogeneous with respect to $\mathcal{A}$ and $r=d+1$, then $P$ is a simplex whose only lattice points are its vertices. 
\end{remark}

\begin{proposition}\label{prop:TFP_factorizes_as_Cayley_sum}
    Let $A \subseteq \Z^{d_1}$ and $B \subseteq \Z^{d_2}$ be homogeneous with respect to $\mathcal{A}$ respectively via affine maps $\pi_1: \Z^{d_1} \to \Z^{\ell}$ and $\pi_2: \Z^{d_2} \to \Z^{\ell}$.
    Then 
    \[
       I_{A \times_{\mathcal{A}} B} = I_{(A_1 \times B_1)*\dots*(A_r \times B_r)}
    \]
    where $A_i$ and $B_i$ are respectively the integer points of the lattice polytopes 
    \[
    P_i := P \cap (\pi_1)_{\R}^{-1}(\mathbf{a}^i),\qquad Q_i := Q \cap (\pi_2)_{\R}^{-1}(\mathbf{a}^i).
    \]
\end{proposition}
\begin{proof}
    Recall that
    \[
    {A \times_{\mathcal{A}}B} = \{ (\mathbf{a}^i_j,\,\mathbf{b}^i_k)\mid i \in [r],\, j \in [s_i],\, k \in [t_i]\}.
    \]
    Let $e_0 := 0$, and let $e_1,\dots,e_{r-1}$ be the standard basis of $\Z^{r-1}$. By the previous lemma, we can assume $A = A_1*\dots*A_r$ and $B=B_1*\dots*B_r$, so the TFP is equal to
    \[
    \bigcup_{i=1}^r A_i \times \{e_{i-1}\} \times B_i \times \{e_{i-1}\} \subseteq \Z^{d_1} \times \Z^{r-1} \times \Z^{d_2} \times \Z^{r-1}.
    \]
    The points $(\mathbf{a},z_1, \mathbf{b},z_2) \in  \Z^{d_1} \times \Z^{r-1} \times \Z^{d_2} \times \Z^{r-1}$ of the TFP all satisfy $z_1 = z_2$, so we can identify them with
    \[
    \bigcup_{i=1}^r A_i  \times B_i \times \{e_{i-1}\} \subseteq \Z^{d_1} \times \Z^{d_2} \times \Z^{r-1}
    \]
    which are the integer points of the polytope
    \[
    (P_1 \times Q_1)*\dots*(P_r \times Q_r). \qedhere
    \]
\end{proof}
Note that by Corollary 2.10 of \cite{Sullivant2007}, one can easily construct a generating set for the toric ideal associated to $(A_1\times B_1)\ast\cdots\ast(A_r\times B_r)$, given homogeneous generating sets of $I_{A}$ and $I_{B}$. 
\begin{example}
    Consider the polytope $\mathcal{P} = (2\Delta_1\times\Delta_1)\ast(\Delta_1\times2\Delta_1)$. This is equivalent to  $T_{1,1,1} \times_{\mathcal{A}} T_{1,1,1}$ where $\mathcal{A}$ is the configuration of vectors $\{(1,0),\,(1,1)\}$ and the multigradings on the lattice points of two copies of $T_{1,1,1}$ are induced by the affine maps
    \[\pi_1(x,y)=(1,y),\qquad \pi_2(x,y) = (1,1-y).\]
  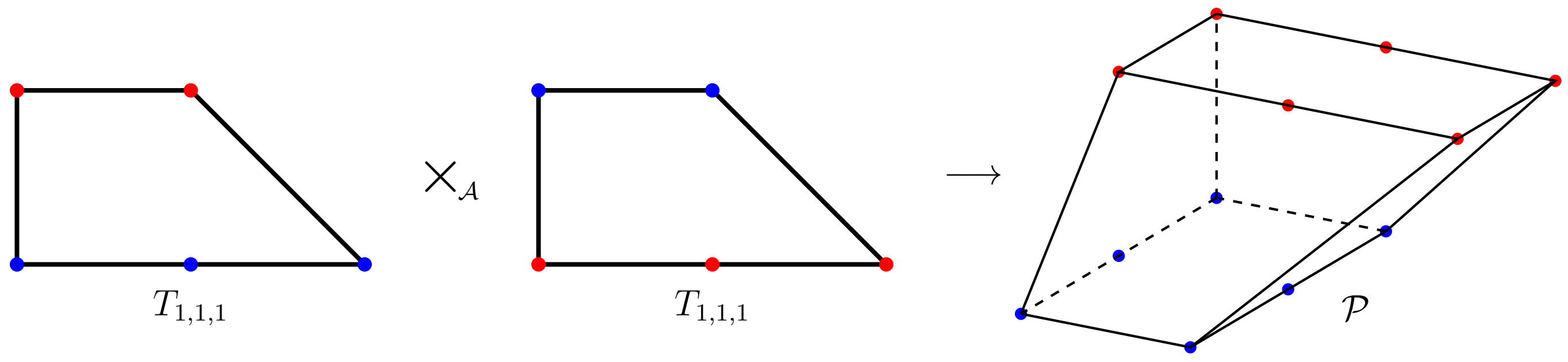
\begin{figure}[H]
    \centering
    \begin{tikzpicture}
    \begin{scope}[line width =1.5pt,scale=0.4]
    \tikzstyle help lines=[color=blue!50,very thin]
    \draw (0,0) -- (5,0) -- (10,0) -- (5,5) -- (0,5) -- (0,0);
    \draw (15,0) -- (25,0) -- (20,5) -- (15,5) -- (15,0);
    \node at (12.5,2.5) {$\bigtimes_{\mathcal{A}}$};
    \node at (27.5, 2.5) {$\longrightarrow$};
    \node at (5,-1.25) {$T_{1,1,1}$} ;
    \node at (20,-1.25) {$T_{1,1,1}$} ;
    \node at (38.5,-1.25) {$\mathcal{P}$} ;
    \filldraw [blue] (0,0) circle (4pt);
    \filldraw [blue] (5,0) circle (4pt);
    \filldraw [blue] (10,0) circle (4pt);
    \filldraw [red] (5,5) circle (4pt);
    \filldraw [red] (0,5) circle (4pt);
    \filldraw [red] (15,0) circle (4pt);
    \filldraw [red] (20,0) circle (4pt);
    \filldraw [red] (25,0) circle (4pt);
    \filldraw [blue] (15,5) circle (4pt);
    \filldraw [blue] (20,5) circle (4pt);
    \end{scope}
    \tdplotsetmaincoords{70}{120} 
    \begin{scope}[tdplot_main_coords, xshift=13.8cm, scale=0.45, yshift=1.7cm]
            
        \coordinate (A) at (0,0,0);
        \filldraw [blue] (A) circle (4pt);        
        \coordinate (Bm) at (5,0,0);
        \filldraw [blue] (Bm) circle (4pt);
        \coordinate (B) at (10,0,0);
        \filldraw [blue] (B) circle (4pt); 
        \coordinate (Cm) at (5,5,0);
        \filldraw [blue] (Cm) circle (4pt);    
        \coordinate (C) at (10,5,0);
        \filldraw [blue] (C) circle (4pt);   
        \coordinate (D) at (0,5,0);
        \filldraw [blue] (D) circle (4pt);    
        \coordinate (A1) at (0,0,5);
        \filldraw [red] (A1) circle (4pt);    
        \coordinate (B1) at (5,0,5);
        \filldraw [red] (B1) circle (4pt);    
        \coordinate (C1m) at (5,5,5);
        \filldraw [red] (C1m) circle (4pt);    
        \coordinate (C1) at (5,10,5);
        \filldraw [red] (C1) circle (4pt);
        \coordinate (D1) at (0,10,5);
        \filldraw [red] (D1) circle (4pt);
        \coordinate (A1m) at (0,5,5);
        \filldraw [red] (A1m) circle (4pt);

        \draw[thick, dashed] (A) -- (B);
        \draw[thick] (B) -- (C) -- (D);
        \draw[thick, dashed] (D) -- (A);

        \draw[thick] (A1) -- (B1) -- (C1) -- (D1) -- (A1);
        
        \draw[thick,dashed] (A) -- (A1);
        \draw[thick] (B) -- (B1);
        \draw[thick] (C) -- (C1);
        \draw[thick] (D) -- (D1);

    \end{scope}
    
\end{tikzpicture}
    \caption{Factorization of $\mathcal{P}$ as a toric fiber product.}
\end{figure}
By definition, the ideal $I_{\mathcal{P}} =I_{T_{1,1,1}} \times_{\mathcal{A}} I_{T_{1,1,1}}$ is the kernel of the ring map
\begin{align*}
            \mathbb{K}[z^1_{00},\dots,z^1_{21},z^2_{00},\dots, z^2_{12}] &\longrightarrow \left( \mathbb{K}[x^1_0,x^1_1,x^1_2,x^2_0,x^2_1]/I_{T_{1,1,1}} \right) \otimes_{\mathbb{K}} \left( \mathbb{K}[y^1_0,y^1_1,y^2_0,y^2_1,y^2_2]/I_{T_{1,1,1}} \right) \\
        z^i_{jk} &\longmapsto x^i_j \otimes y^i_k.
\end{align*}
We can calculate that it is generated by the following 32 quadrics
\begin{align*}
& (z^2_{11})^2 - z^2_{10}z^2_{12},\;z^1_{21}z^2_{02} - z^1_{11}z^2_{12},\;z^1_{20}z^2_{00} - z^1_{10}z^2_{10} 
,\; \\ &z^2_{02}z^2_{11} - z^2_{01}z^2_{12},\;z^1_{20}z^2_{02} - z^1_{10}z^2_{12},\;z^1_{11}z^2_{00} - z^1_{00}z^2_{11} 
,\; \\ &z^2_{01}z^2_{11} - z^2_{00}z^2_{12},\;z^1_{11}z^2_{02} - z^1_{01}z^2_{12},\;z^1_{10}z^2_{00} - z^1_{00}z^2_{10} 
,\; \\ &z^1_{21}z^2_{11} - z^1_{20}z^2_{12},\;z^1_{10}z^2_{02} - z^1_{00}z^2_{12},\;z^1_{01}z^2_{00} - z^1_{00}z^2_{01} 
,\; \\ &z^1_{11}z^2_{11} - z^1_{10}z^2_{12},\;(z^2_{01})^2 - z^2_{00}z^2_{02},\;z^1_{11}z^1_{20} - z^1_{10}z^1_{21} 
,\; \\ &z^1_{01}z^2_{11} - z^1_{00}z^2_{12},\;z^1_{21}z^2_{01} - z^1_{10}z^2_{12},\;z^1_{01}z^1_{20} - z^1_{00}z^1_{21} 
,\; \\ &z^2_{02}z^2_{10} - z^2_{00}z^2_{12},\;z^1_{20}z^2_{01} - z^1_{10}z^2_{11},\;(z^1_{11})^2 - z^1_{01}z^1_{21} 
,\; \\ &z^2_{01}z^2_{10} - z^2_{00}z^2_{11},\;z^1_{11}z^2_{01} - z^1_{00}z^2_{12},\;z^1_{10}z^1_{11} - z^1_{00}z^1_{21} 
,\; \\ &z^1_{21}z^2_{10} - z^1_{20}z^2_{11},\;z^1_{10}z^2_{01} - z^1_{00}z^2_{11},\;(z^1_{10})^2 - z^1_{00}z^1_{20} 
,\; \\ &z^1_{11}z^2_{10} - z^1_{10}z^2_{11},\;z^1_{01}z^2_{01} - z^1_{00}z^2_{02},\;z^1_{01}z^1_{10} - z^1_{00}z^1_{11} 
,\; \\ &z^1_{01}z^2_{10} - z^1_{00}z^2_{11},\;z^1_{21}z^2_{00} - z^1_{10}z^2_{11}
\in \mathbb{K}[z^1_{00},\dots,z^2_{12}].    
\end{align*}
\end{example}

\subsection{Classification of 3D codimension zero TFPs with ML degree one}
Using results from the previous subsection, we now classify all $3D$ polytopes with ML degree one that are obtained as codimension zero TFPs of lower-dimensional polytopes. First, we recall that \cite{vonBothmer2010} proved that the only polytopes of dimension $\leq 2$ with ML degree one are equivalent to a member of the family
\[
\left\{ T_{a,b,d} := \mathrm{conv}\begin{pmatrix}
    0 & a+db & 0 & a  \\
    0 &  0   & b & b 
\end{pmatrix}:\ (a,\,b,\,d) \in \mathbb{N}^3\right\} .
\]
By compiling all possible ways in which the polytopes of the form $T_{a,b,d}$ can be homogeneously graded with respect to some linearly independent $\mathcal{A}$, we can obtain our classification in Theorem~\ref{thm:classification}. Note that the classification yields a complete description of $3D$ codimension zero TFPs of lower-dimensional polytopes, since they satisfy the following property.

\begin{proposition}\label{prop:3D_TFP}
Let $P$ be a $3$-dimensional polytope that is a codimension zero TFP of lower dimensional polytopes. Then, one of the following is true:
\begin{enumerate}
    \item $P$ is the Segre product of a polygon and an interval;
    \item $P$ has ML degree $1$. 
\end{enumerate}
\end{proposition}
\begin{proof} Take a lattice polytope $P \subseteq  \mathbb{R}^2$ whose integer points are homogeneous with respect to $\mathcal{A}$. If $r+1=3$, then $P=\Delta_2$ by Remark~\ref{remark:simplex} and it cannot be a factor of any $3$-dimensional TFP multigraded by $\mathcal{A}$. If $r+1=2$, then the integer points lie on two distinct lines, so up to a unimodular transformation the polytope is $P = T_{a,1,d}$ for some $a,\,d \in \mathbb{N}$; hence, every TFP with $r+1=2$ has ML degree one. If $r+1=1$, then it is a Segre product. 
\end{proof}

\begin{theorem}\label{thm:classification}
Let $L$ be a $3$-dimensional polytope obtained as a codimension zero TFP of lower-dimensional polytopes with ML degree one, and let $c$ be any vector of positive scalings associated with the integer points of $L$. Then, there exists a polytope $P$ in the family
    \begin{align*}
    \mathcal{F} := \; & \{ \left( a\Delta_1 \times b\Delta_1\right) * \left( a'\Delta_1 \times b'\Delta_1    \right) \mid a,\,b,\,a',\,b' \in \mathbb{N} \} \\
    & \cup ~ \{T_{a,b,d}\times h\Delta_1 \mid a,b,d,h \in \mathbb{N} \}
    \end{align*}
    such that
    \[
    X_{L,c} = X_{P,c}.
    \]
\end{theorem}
\begin{proof}
    Since $L$ is $3$-dimensional, at least one of its factors in the TFP is $2$-dimensional. By the classification of polygons with ML degree one, we may assume that this factor is $T_{a,b,d}$, for some $a,b,d \in \mathbb{N}$. By Lemma \ref{lemma:CayleySum}, a multigrading with respect to some linearly independent vectors $\mathcal{A} = \{\mathbf{a}^1,\dots,\mathbf{a}^r\}$ induces a partition of its integer points in $r$ parallel subspaces of codimension $r-1$. 
    
    If $r=3$, by Remark \ref{remark:simplex} the polytope $T_{a,b,d}$ is unimodularly equivalent to $\Delta_2$. Since the second factor of the TFP is also homogeneous with respect to $\mathcal{A}$, which is a configuration of three linearly independent vectors, it cannot have dimension strictly less than $2$. Therefore, the second factor is also unimodularly equivalent to $\Delta_2$. Since
    \[
    \Delta_2 \times_{\mathcal{A}} \Delta_2 = \Delta_2
    \]
    is not a $3$-dimensional polytope, this case is impossible.

    If $r=2$, then the multigrading induces a partition of the integer points of $T_{a,b,d}$ into two parts lying on two parallel lines. 
    Consider the case where $b > 1$. The edge of $T_{a,b,d}$ whose vertices are $(0,0)$ and $(0,b)$ contains at least $3$ integer points. Since $r=2$, at least two of these points must lie on the same line, hence the integer points of $T_{a,b,d}$ are contained in two vertical lines. Up to a unimodular transformation, these are the integer points of the polytope $T_{b,1,0}$. 
    Consider now the case where $b=1$. If $a+d \geq 2$, the previous argument shows that the parallel lines must be horizontal. If $a+d=1$, the polygon is either $\Delta_2$ or the unit square.
    Up to an affine unimodular transformation, every admissible grading for $r=2$ induces the partition
    \[
   \Big\{(a+d)\Delta_1\times\{0\},\,a\Delta_1\times\{1\} \Big\}.
    \]
    Hence, for $r=2$, all the $3$-dimensional TFPs have the form
    \begin{align*}
    T_{a,1,d} \times_{\mathcal{A}_1} T_{a',1,d'} &= \Big( (a+d)\Delta_1 \times (a'+d')\Delta_1 \Big)* \Big( a\Delta_1 \times a'\Delta_1 \Big), \text{ or}\\
    T_{a,1,d} \times_{\mathcal{A}_2} T_{a',1,d'} &= \Big( (a+d)\Delta_1 \times a'\Delta_1 \Big) * \Big( a\Delta_1 \times (a'+d')\Delta_1 \Big)
    \end{align*}
    for some sets $\mathcal{A}_1,\,\mathcal{A}_2 \subseteq \Z^{\ell}$ made of two linearly independent vectors.
    Up to equivalence, by varying the values of $a,\,a',d,\,d' \in \mathbb{N}$ we get all the $3$-dimensional polytopes of the form 
    \[
    \left( a\Delta_1 \times b\Delta_1\right) * \left( a'\Delta_1 \times b'\Delta_1    \right) 
    \]
    for $a,\,a',\,b,\,b' \in \mathbb{N}$.

    If $r=1$, the TFP is just the Segre product. Hence, we obtain all the polytopes of the form
    \[
    T_{a,b,d}\times h\Delta_1 
    \]
    for $a,\,b,\,d,\,h \in \mathbb{N}$.
\end{proof}

\begin{remark}
Not all elements of the family $\mathcal{F}$ in Theorem~\ref{thm:classification} are $3$-dimensional polytopes. In particular, if we remove all of them from $\mathcal{F}$, there would remain the collection of ML degree one polytopes of dimensions less than or equal to $2$, up to unimodular transformation.
\end{remark}
\begin{example}
    Consider the polytope $\mathcal{P} \subseteq \R^3$, constructed as
    \[
    \operatorname{conv}  
    \begin{pmatrix}
        0 & 0 & 1 & 1 & 0 & 1 & 2\\
        0 & 1 & 0 & 1 & 1 & 1 & 1 \\ 
        0 & 0 & 0 & 0 & 1 & 1 & 1\\
    \end{pmatrix}.
    \]
    We notice that its lattice points lie in the planes $z=0$ and $z=1$.
      \begin{figure}[H]
    \centering
    \begin{tikzpicture}
    \tdplotsetmaincoords{70}{120} 
    \begin{scope}[tdplot_main_coords, xshift=11cm, scale=0.45, yshift=1.25cm]
            
        \coordinate (A) at (0,0,0);
        \filldraw [blue] (A) circle (4pt);        
        \coordinate (B) at (5,0,0);
        \filldraw [blue] (B) circle (4pt);    
        \coordinate (C) at (5,5,0);
        \filldraw [blue] (C) circle (4pt);    
        \coordinate (D) at (0,5,0);
        \filldraw [blue] (D) circle (4pt);    
        \coordinate (A1) at (0,5,5);
        \filldraw [red] (A1) circle (4pt);    
        \coordinate (B1) at (5,5,5);
        \filldraw [red] (B1) circle (4pt);    
        \coordinate (C1) at (10,5,5);
        \filldraw [red] (C1) circle (4pt);

        \draw[thick, dashed] (A) -- (B);
        \draw[thick] (B) -- (C) -- (D);
        \draw[thick, dashed] (D) -- (A);

        \draw[thick] (A1) -- (B1) -- (C1);
        
        \draw[thick,dashed] (A) -- (A1);
        \draw[thick] (D) -- (A1);
        \draw[thick] (B) -- (C1);
        \draw[thick] (C) -- (C1);

    \end{scope}
    
\end{tikzpicture}
    \caption{The lattice points of $\mathcal{P}$; those contained in $z=0$ are colored in blue, those contained in $z=1$ are colored in red.}
\end{figure}
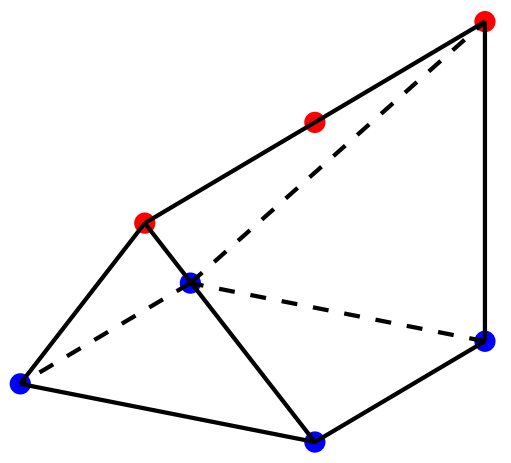
As a consequence, $\mathcal{P}$ is the Cayley sum of a square $\Delta_1^2$ and of a horizontal line $2\Delta_1 \times \{1\}$. It corresponds to the eighth item of Table \ref{tb:polyhedra} with $a=a'=d'=1$, therefore there exists a positive scaling for which the associated embedded scaled toric variety has ML degree $1$. Indeed, by Proposition \ref{prop:TFP_factorizes_as_Cayley_sum}, we have that $\mathcal{P}$ is equivalent to
\begin{align*}
(\Delta_1 \times \Delta_1) \ast (\{0\}\times 2\Delta_1)
&\cong \Delta_2 \times_{\mathcal{A}} T_{1,1,1},
\end{align*}
where $\mathcal{A}$ is the configuration of vectors $\{(1,0),\,(1,1)\}$ and the multigradings on the lattice points of $\Delta_2$ and $T_{1,1,1}$ are respectively induced by the affine maps
    \[\pi_1(x,y)=(1,y),\qquad \pi_2(x,y) = (1,1-y).\]

  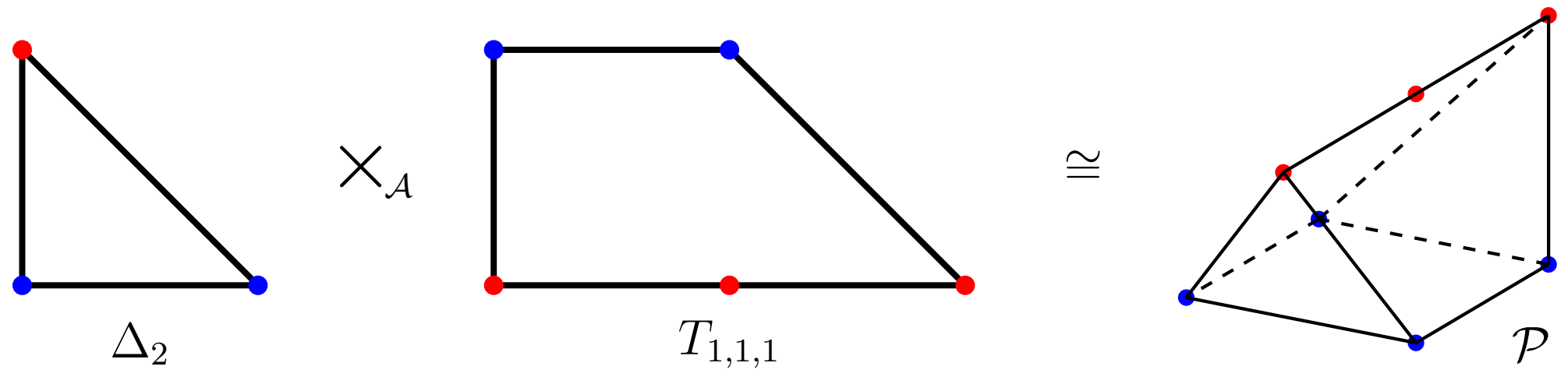
\begin{figure}[H]
    \centering
    \begin{tikzpicture}
    \begin{scope}[line width =1.5pt,scale=0.4]
    \tikzstyle help lines=[color=blue!50,very thin]
    \draw (0,0) -- (5,0) -- (0,5) -- (0,0);
    \draw (10,0) -- (20,0) -- (15,5) -- (10,5) -- (10,0);
    \node at (7.5,2.5) {$\bigtimes_{\mathcal{A}}$};
    \node at (22.5, 2.5) {$\cong$};
    \node at (2.5,-1.25) {$\Delta_2$} ;
    \node at (15,-1.25) {$T_{1,1,1}$} ;
    \node at (32,-1.25) {$\mathcal{P}$} ;

    \filldraw [blue] (0,0) circle (4pt);
    \filldraw [blue] (5,0) circle (4pt);
    \filldraw [red] (0,5) circle (4pt);
    \filldraw [red] (10,0) circle (4pt);
    \filldraw [red] (15,0) circle (4pt);
    \filldraw [red] (20,0) circle (4pt);
    \filldraw [blue] (10,5) circle (4pt);
    \filldraw [blue] (15,5) circle (4pt);
    \end{scope}
    \tdplotsetmaincoords{70}{120} 
    \begin{scope}[tdplot_main_coords, xshift=11cm, scale=0.45, yshift=1.25cm]
            
        \coordinate (A) at (0,0,0);
        \filldraw [blue] (A) circle (4pt);        
        \coordinate (B) at (5,0,0);
        \filldraw [blue] (B) circle (4pt);    
        \coordinate (C) at (5,5,0);
        \filldraw [blue] (C) circle (4pt);    
        \coordinate (D) at (0,5,0);
        \filldraw [blue] (D) circle (4pt);    
        \coordinate (A1) at (0,5,5);
        \filldraw [red] (A1) circle (4pt);    
        \coordinate (B1) at (5,5,5);
        \filldraw [red] (B1) circle (4pt);    
        \coordinate (C1) at (10,5,5);
        \filldraw [red] (C1) circle (4pt);

        \draw[thick, dashed] (A) -- (B);
        \draw[thick] (B) -- (C) -- (D);
        \draw[thick, dashed] (D) -- (A);

        \draw[thick] (A1) -- (B1) -- (C1);
        
        \draw[thick,dashed] (A) -- (A1);
        \draw[thick] (D) -- (A1);
        \draw[thick] (B) -- (C1);
        \draw[thick] (C) -- (C1);

    \end{scope}

\end{tikzpicture}
    \caption{Factorization of $\mathcal{P}$ as a toric fiber product.}
\end{figure}
Furthermore, we may calculate its associated ideal  $I_{\mathcal{P}} = I_{\Delta_2} \times_{\mathcal{A}} I_{T_{1,1,1}}$ as the kernel of the ring map
\begin{align*}
            \mathbb{K}[z^1_{00},\dots,z^1_{11},z^2_{00},\dots, z^2_{02}] &\longrightarrow \left( \mathbb{K}[x^1_0,x^1_1,x^2_0]/I_{\Delta_2} \right) \otimes_{\mathbb{K}} \left( \mathbb{K}[y^1_0,y^1_1,y^2_0,y^2_1,y^2_2]/I_{T_{1,1,1}} \right) \\
        z^i_{jk} &\longmapsto x^i_j \otimes y^i_k
\end{align*}
where 
\begin{align*}
I_{\Delta_2} &= (0) \subseteq \mathbb{K}[x_0^1,x_1^1,x_0^2], \\
I_{T_{1,1,1}} &= \left( (y^2_1)^2 - y^2_0 y^2_2,\, y^1_1 y^2_1 - y^1_0 y^2_2,\, y^1_1 y^2_0 - y^1_0 y^2_1 \right) \subseteq \mathbb{K}[y_0^1,\dots, y_2^2].
\end{align*}
In particular, $I_{\mathcal{P}}$ is generated by the following 6 quadrics
\begin{align*}
    & (z^2_{01})^2 - z^2_{00}z^2_{02},\; 
    z^1_{11}z^2_{01} - z^1_{10}z^2_{02},\; 
    z^1_{01}z^2_{01} - z^1_{00}z^2_{02},\; \\
    &z^1_{11}z^2_{00} - z^1_{10}z^2_{01},\; 
    z^1_{01}z^2_{00} - z^1_{00}z^2_{01},\; 
    z^1_{01}z^1_{10} - z^1_{00}z^1_{11}
    \in \mathbb{K}[z^1_{00},\dots,z^1_{11},z^2_{00},\dots,z^2_{02}].
\end{align*}
\end{example}

\section{3D polytopes with ML degree higher than one but with ML degree one faces}
In this section we consider one family of 3D polytopes with ML degree higher than one but with ML degree one faces. The computation of the ML degree of this family is based on an explicit Euler characteristic calculation that reduces the problem to counting roots of certain univariate polynomials in a special way, this is presented in section \ref{sec:reduction}. This argument is similar to the
one by Telen and Wiesmann \cite{telen2024} for Hirzebruch surfaces. These results allow us to find different scalings that yield
a complete list of the ML degrees that are attained.

\subsection{Euler characteristic computation by counting roots of univariate polynomials}
By Theorem \ref{thm:HuhSturmfels} the ML degree of a toric variety equals the (signed) Euler characteristic of a hypersurface complement in the torus $(\C^{\ast})^d$. As $\chi ((\C^\ast)^d) = 0$ we can equivalently compute the Euler characteristic of the hypersurface. Consider a polynomial $f\in \C[x_1,\dots,x_{d-1},t]$. We can regard $f$ as a parametric univariate polynomial in $t$ with coefficients in $\C[x_1,\dots,x_{d-1}]$, written as
\[
    f= q_0(x)+q_1(x)t+\dots +q_n(x)t^n
\]
where $q_{i}(x)\in \C[x_1,\ldots,x_{d-1}]$.
Specialization of coefficients defines a map 
\[
    \phi_f \colon (\C^\ast)^{d-1} \to \C[t] ,
\]
which takes a choice of coefficients $c\in (\C^{\ast})^{d-1}$ to the univariate
polynomial $f_c(t)$ with complex coefficients $q_{i}(c)$.
For $k\in \Z_{\geq0}$ define $O_k$ to be the set of polynomials in $\C[t]$ with exactly $k$ distinct nonzero roots. We write an expression for the Euler characteristic of a hypersurface in the following lemma. 

\begin{lemma}       
    For a polynomial $f\in \C[x_1,\dots,x_{d-1},t]$ the Euler characteristic can be computed by
    \[
        \chi (V(f)) = \sum _{i=1}^{n} i \cdot \chi(\phi_f^{-1}(O_i))
    \]
\end{lemma}
\begin{proof}
    Let $\pr_x \colon (\C^\ast)^d \to (\C^\ast)^{d-1}$ be the projection to the $x$-coordinates. Restricting $\pr_x$ to $V(f)\cap\pr_x^{-1}(\phi_f^{-1}(O_l))$ gives a fibration with discrete fiber of cardinality $l$ and base $\phi_f^{-1}(O_l)$. By multiplicativity of the Euler characteristic on fibrations $\chi(V(f)\cap\pr_x^{-1}(\phi_f^{-1}(O_l)))= l*\chi(\phi_f^{-1}(O_l))$. The formula follows from additivity of the Euler characteristic. Note that the fiber over $\phi_f^{-1}(0)$ is $\C^\ast$, thus it does not appear in the summation.
\end{proof}
Note that $\bigcup_i \phi_f^{-1}(O_i) \cup  \phi_f^{-1}(0)$ induces a partition of $(\C^\ast)^{d-1}$ into locally closed sets. They can be defined precisely using subdiscriminants which can be obtained from the classical Sylvester formula \cite{chen2025}.
In general one cannot expect the previous lemma to carry much weight in analyzing ML degrees of hypersurfaces associated to scaled projective toric varieties. However, we want to apply it to concrete examples where $f$ specializes to univariate polynomials with small degree.
\begin{corollary}\label{degree 1 case}
    Assume $f= q_0(x)+q_1(x)t\in\C[x_1,\dots,x_{d-1},t]$, then 
    \[
        \chi(V(f)) = -\chi(V(q_0(x)q_1(x)))
    \]
\end{corollary}
\begin{proof}
    As $f$ specializes to a degree one polynomial we see that $\phi_f^{-1}(O_1)= (\C^\ast)^{d-1}\setminus V(q_0(x)q_1(x))$.
\end{proof}
This result is the case of one variable of a theorem of Fevola and Matsubara-Heo
\cite[Theorem~2.2]{fevola2024}, proved in the study of Euler discriminants.
It has been applied to the study of ML degrees in \cite{hosten2026}, where the
Euler characteristic for $\mathbb{P}^1 \times \mathbb{P}^1\times \mathbb{P}^n$ is
reduced to that of a union of $n+1$ quadric curves in $(\C^\ast)^2$. A similar
fibration argument appears in \cite{clarke2024} for independence models
$\mathbb{P}^{m}\times\mathbb{P}^{n}$, whose ML degrees are Euler characteristics
of hyperplane arrangement complements. Corollary~\ref{cor:deg2case} is proved by
an argument similar to the previous one, applied to a product of two polynomials
that are linear in~$t$.
\begin{corollary} \label{cor:deg2case}
    Assume $f= (q_0(x)+q_1(x)t)(q_2(x)+q_3(x)t)\in \C[x,t]$, set $D_f=q_0q_3-q_1q_2 \in \C[x]$ and $M_1 = \{x\in \C^\ast \,: \text{exactly one of } q_0,q_1,q_2,q_3,D_f \text{ vanishes at } x \}$, $M_2=\C^\ast \setminus V(q_0q_1q_2q_3D_f)$, then
    \[
    \chi(V(f)) = 2\chi(M_2)+\chi(M_1)
    \]
\end{corollary}
\begin{proof}
    The preimage  of $O_2$ under $\phi_f$ is given by $\C^\ast \setminus V(q_0q_1q_2q_3D_f)$. It is easy to see that the fiber of $O_1$ is given by $M_1$.
\end{proof}
The Newton polytope $P_f$ associated to the polynomial $f\in\C[x_1,\ldots,x_{d-1},t]$ is the convex hull of the
set of exponent vectors with respect to the monomials
in the variables $x_1,\ldots,x_{d-1},t$ that appear in $f$.
We collect these exponent vectors as the columns of a matrix $A$ and
require that the set of lattice points in $P_f$ coincides
with the set of columns in $A$. From $f$ we also get a vector
$c$ of scalings associated to each column in $A$. The pair
$(A,c)$ determines a scaled projective toric variety $X_{(A,c)}$.
In this section we also denote it by $X_{f}$ since it is readily
obtained from $f$.

Using the previous results we get another direct proof for the ML degree of rational normal scrolls.
\begin{lemma}\cite{amendola2019}
    For the polytope defined by $G = g_0(x)+g_1(x)t_1+g_2(x)t_2+\dots g_n(x)t_n \in \C[x,t_1,\dots,t_n]$. The ML degree for the variety $X_{G}$ is given by the number of distinct nonzero roots of the polynomials $g_i(x)$.
\end{lemma}
\begin{proof}
     By Theorem \ref{thm:HuhSturmfels}, the ML degree of $X_{G}$ is $(-1)^{n+1}\chi((\C^\ast)^{n+1}\setminus V(G))$. Applying Corollary~\ref{degree 1 case} for each $t_i$ consecutively yields $\mldeg X_{G} = \chi(V(g_0(x)g_1(x)\cdots g_n(x)))$. 
\end{proof}

\subsection{ML degrees of 3D polytopes with vertices in $\Delta_1\times \Delta_1\times \Z_{\geq0}$} \label{sec:reduction}
We  employ the calculations in the previous subsection to compute ML degrees of certain 3D polytopes. To make computations possible, we 
restrict to the case of 3D lattice polytopes
defined as the convex hull
of a set of lattice points in $\Delta_1\times \Delta_1\times \Z_{\geq 0}$.
Consider the 3D polytope $P_{f}$ defined as the Newton polytope of the polynomial
\[
    f = f_0(x) + f_1(x)s + f_2(x)t + f_3(x)st
\] 
where $f_i(x)=\sum_{j\geq 0}c_{ij}x^j\in \C[x]$.  Let $b_i=\deg(f_i)$, $a_i=\ord(f_i)$ and $d_i=b_i-a_i$ for all $i\in \{0,1,2,3\}$. Recall that the order
of a univariate polynomial is the multiplicity of zero as a root. Then each $f_i$ is equal to $x^{a_i}$ times a polynomial $g_i$ of degree $d_i$, where $x$ does not divide $g_i$. We compute the ML degree of $X_{f}$ as the signed Euler characteristic of a hypersurface complement defined by $V(f)$. Using the results of the previous section yields
\begin{align}\label{computation 3d polytopes}
    \mldeg X_{f} &=  (-1)^3\chi((\C^\ast)^3\setminus V(f)) \nonumber \\
      &\overset{\ref{degree 1 case}}{=}  -\chi ( V((f_0(x) + f_1(x)s)(f_2(x)+f_3(x)s))) \nonumber \\
    &\overset{\ref{cor:deg2case}}{=} -2\chi(M_2)-\chi(M_1)
\end{align}
with the sets $M_1,M_2$ defined as in Corollary \ref{cor:deg2case}. 

So far we have not made assumptions about  $f_i$'s to be nonzero. If $f_0=0$, then $M_2=\emptyset $ and $M_1 = \C^\ast \setminus V(f_1f_2f_3)$, the 
ML degree is then given by the number of nonzero roots among the $f_i$.  From now on we assume all the $f_i$ are nonzero. The $f_i$
can be thought of as posts standing on the points $(0,0,0),(1,0,0),(0,1,0),(1,1,0)$  as shown in Figure~\ref{fig:Pf-general}. The $d_i$ are the lengths
of the posts of $P_f$. 
\begin{definition}\label{def:tubular}
We call $P_f$ \textit{tubular} if $f_0,f_1,f_2,f_3$ are nonzero and neither $f_0f_3$ nor $f_1f_2$ is a monomial, that is, $d_0+d_3\geq 1$ and $d_1+d_2\geq 1$.
\end{definition}
Tubular polytopes have eight vertices $w_{000},w_{001},\dots,w_{111}$ given by the coordinates $w_{ab0}= (a,b,\ord (f_{a+2b}))$ and $w_{ab1}=(a,b,\deg (f_{a+2b}))$ for $a,b \in \{0,1\}$, whenever all $d_i$ are positive. 
\begin{figure}[ht]
  \centering
  \definecolor{postgreen}{RGB}{0,190,40}
\subfigure[A general tubular polytope $P_f$ inside the tube $\Delta_1\times\Delta_1\times\Z$.]{%
\begin{tikzpicture}[
    scale=1.25,
    x={(1.1cm,-0.75cm)}, y={(1.75cm,0.6cm)}, z={(0cm,1.75cm)},
    edge/.style   = {line width=1.3pt, line cap=round, line join=round},
    hidden/.style = {edge, dash pattern=on 4pt off 3pt},
    post/.style   = {edge, postgreen, line width=1.8pt},
    hiddenpost/.style = {post, dash pattern=on 4pt off 3pt},
    tube/.style   = {line width=1.1pt, line cap=round, dash pattern=on 0pt off 3.5pt},
    vertex/.style = {circle, fill=black, inner sep=0pt, minimum size=6pt},
    lbl/.style    = {font=\normalsize, inner sep=1pt},
    flbl/.style   = {lbl, text=postgreen}
  ]
  \coordinate (w000) at (0,0,0.2);  \coordinate (w001) at (0,0,1.55);
  \coordinate (w100) at (1,0,0.0);  \coordinate (w101) at (1,0,0.7);
  \coordinate (w010) at (0,1,0.4);  \coordinate (w011) at (0,1,1.4);
  \coordinate (w110) at (1,1,0.2);  \coordinate (w111) at (1,1,0.55);
  \foreach \a/\b/\lo/\hi in {0/0/0.2/1.55, 1/0/0.0/0.7, 0/1/0.4/1.4, 1/1/0.2/0.55}{
    \draw[tube] (\a,\b,\hi)  -- (\a,\b,2.1);
  }
  \foreach \a/\b/\lo in {0/0/0.2, 1/0/0.0, 1/1/0.2}
    \draw[tube] (\a,\b,-0.5) -- (\a,\b,\lo);
  \draw[hidden] (w000) -- (w010) -- (w110);
  \draw[edge] (w000) -- (w100) -- (w110);
  \draw[edge] (w001) -- (w011) -- (w111) -- (w101) -- cycle;
  \draw[hiddenpost] (w010) -- (w011);
  \draw[post] (w000) -- (w001);
  \draw[post] (w100) -- (w101);
  \draw[post] (w110) -- (w111);
  \foreach \v in {w000,w001,w010,w011,w100,w101,w110,w111}
    \node[vertex] at (\v) {};
  \node[lbl, left=5pt]        at (w000) {$w_{000}$};
  \node[lbl, left=5pt]        at (w001) {$w_{001}$};
  \node[lbl, left=5pt]        at (w100) {$w_{100}$};
  \node[lbl, right=6pt, yshift=-5pt] at (w101) {$w_{101}$};
  \node[lbl, right=6pt, yshift=5pt]  at (w010) {$w_{010}$};
  \node[lbl, above left=2pt]  at (w011) {$w_{011}$};
  \node[lbl, below right=2pt] at (w110) {$w_{110}$};
  \node[lbl, above right=2pt] at (w111) {$w_{111}$};
  \node[flbl, left=4pt]  at (0,0,0.875) {$f_0$};
  \node[flbl, right=4pt] at (1,0,0.35)  {$f_1$};
  \node[flbl, left=4pt]  at (0,1,0.9)   {$f_2$};
  \node[flbl, right=4pt] at (1,1,0.375) {$f_3$};
\end{tikzpicture}%
\label{fig:Pf-general}}\hspace{2cm}
\subfigure[The polytope $P_f$ of Example~\ref{ex:simplex-cut}; the post $f_1$ is the single point $c_2$.]{%
\begin{tikzpicture}[
    scale=1.25,
    x={(1.1cm,-0.75cm)}, y={(1.75cm,0.6cm)}, z={(0cm,1.75cm)},
    edge/.style   = {line width=1.3pt, line cap=round, line join=round},
    hidden/.style = {edge, dash pattern=on 4pt off 3pt},
    post/.style   = {edge, postgreen, line width=1.8pt},
    hiddenpost/.style = {post, dash pattern=on 4pt off 3pt},
    vertex/.style = {circle, fill=black, inner sep=0pt, minimum size=6pt},
    lbl/.style    = {font=\normalsize, inner sep=1pt},
    flbl/.style   = {lbl, text=postgreen}
  ]
  \coordinate (C1) at (0,0,0);
  \coordinate (C2) at (1,0,0);
  \coordinate (C3) at (0,1,0);
  \coordinate (C4) at (1,1,0);
  \coordinate (C5) at (0,0,1);
  \coordinate (C6) at (0,1,1);
  \coordinate (C7) at (1,1,1);

  \draw[hiddenpost] (C3) -- (C6);
  \draw[post] (C1) -- (C5);
  \draw[post] (C4) -- (C7);
  
  \draw[edge] (C1) -- (C2) -- (C4);
  \draw[edge] (C5) -- (C6) -- (C7) -- cycle;
  \draw[edge] (C5) -- (C2) -- (C7);
  \draw[hidden] (C1) -- (C3) -- (C4);

  \foreach \v in {C1,C3,C4,C5,C6,C7}
    \node[vertex] at (\v) {};
  \node[vertex, fill=postgreen] at (C2) {};
  \node[lbl, left=5pt]        at (C1) {$c_1$};
  \node[lbl, below=5pt]       at (C2) {$c_2$};
  \node[lbl, above left=2pt]  at (C3) {$c_3$};
  \node[lbl, right=5pt]       at (C4) {$c_4$};
  \node[lbl, left=5pt]        at (C5) {$c_5$};
  \node[lbl, above=5pt]       at (C6) {$c_6$};
  \node[lbl, right=5pt]       at (C7) {$c_7$};
  \node[flbl, left=4pt]  at (0,0,0.5) {$f_0$};
  \node[flbl, right=6pt, yshift=-7pt] at (C2) {$f_1$};
  \node[flbl, left=4pt]  at (0,1,0.35) {$f_2$};
  \node[flbl, right=4pt] at (1,1,0.5) {$f_3$};
\end{tikzpicture}%
\label{fig:Pf-example}}
  \caption{Tubular polytopes. The posts $f_0,\dots,f_3$ (in green) are the vertical edges over the corners of the square $\Delta_1\times\Delta_1$; their lengths are $d_0,\dots,d_3$.}\label{fig:Pf}
\end{figure}
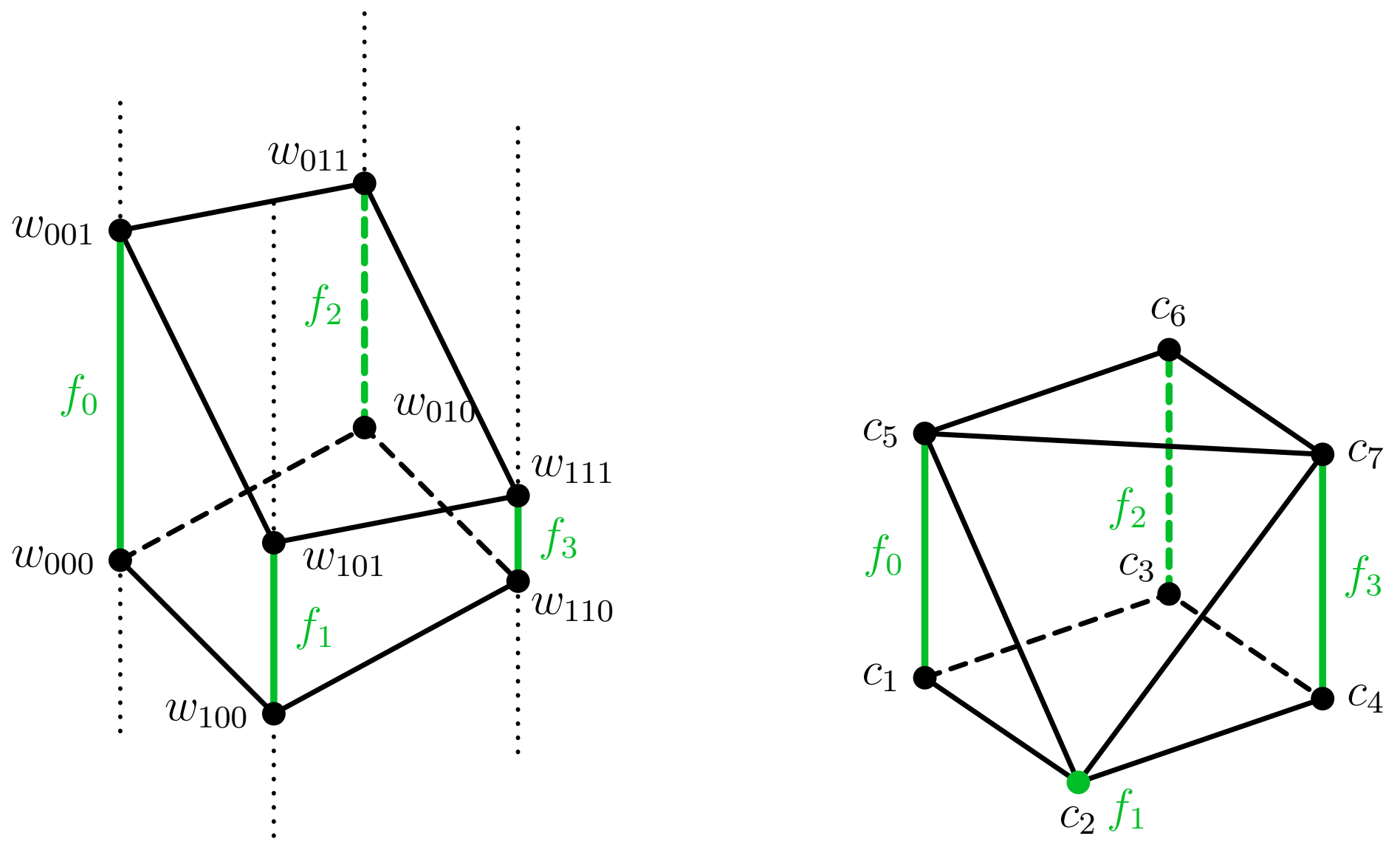
\begin{example} \label{ex:simplex-cut}
Consider the polytope $P_f$, shown in Figure~\ref{fig:Pf-example}, where $f_0=c_1+c_5x$, $f_1=c_2$, $f_2=c_3+c_6x$,
$f_3=c_4+c_7x$ and 
\[f= f_0+ f_1 s + f_2t+ f_3st.\]
The choice of coefficients $c_i=1$ for all $i\in\{1,\ldots,7\}$ yields 
$f_0=f_2=f_3=1+x$, $f_1=1$ and $D_f=x(1+x)$. On the one hand, we see that $M_1$ 
is the empty set because the $f_i$, $i\neq 1$, are identical of degree one, so 
there are no points in which just one of the polynomials in $f_0,\ldots,f_3,D_f$ vanishes. On the other hand, $M_2= \C^{\ast}\setminus V((1+x)^4x)=\C^\ast\setminus \{-1\}$. Hence $\mldeg(X_{f})= -2 \chi(M_2) - \chi(M_1) = -2(-1) - 0 = 2.$ Note that the restrictions
to the faces of this polytope yield polygons with scalings that give ML degree one for all the faces; nevertheless, the full polytope
has ML degree two.
If the coefficients $c_i$ take different values, for instance, $c_1=c_3$, $c_2=c_4$, $c_5=c_6$ and $c_1c_7\neq c_2 c_5$, we obtain
\[
    \begin{split}
        f_0 &= f_2,\\
        D_f &= f_0(f_3-f_1) = f_0c_7x,\\
        M_2 &= \C^* \setminus \{-c_1/c_5,-c_2/c_7\},\\
        M_1 &= \{-c_2/c_7\}.
    \end{split}
\]
The ML degree in this case is
$   \mldeg X_f = -2(-2) - 1 = 3$.
\end{example}
\begin{definition}\label{def:squarebase}
   We say that $P_f$ has a \textit{square base} if exactly one of the sets $\{w_{000},w_{010},w_{100},w_{110}\}$ and  $\{ w_{001},w_{011},w_{101},w_{111}\}$ is a facet of the polytope $P_f$. If both are facets, we say that $P_f$ is \textit{square}. 
\end{definition}
 The polytope in Example~\ref{ex:simplex-cut} has a square base.
 \begin{lemma}\label{lem:faces}
Let $P_f$ be a tubular polytope, then every face of $P_f$ admits a scaling with ML degree one.
\end{lemma}
\begin{proof}
Let $F$ be a facet of $P_f$, then $F$ is contained in a hyperplane. There are only two possibilities for $F$, a triangle or a trapezoid (possibly a parallelogram). In each case we can find a unimodular transformation that sends the face to one of the polygons
in the list of polygons with ML degree one, see Figure~\ref{fig:2Dmldegree-one}.
Since the facets admit a scaling with ML degree one, so does every face of $P_f$
\cite{coons2021,amendola2025}.
\end{proof}
\begin{remark}
    The condition for a square base implies exactly one of the following equalities holds $\ord(f_0)+\ord(f_3) = \ord(f_1)+\ord(f_2)$ or $\deg (f_0) + \deg (f_3) = \deg (f_1) + \deg (f_2)$. By using the symmetry of the polytope we can always assume that a polytope satisfying the square base condition contains the facet $\{w_{000},w_{010},w_{100},w_{110} \}$. 
\end{remark}
\begin{lemma} \label{nonvanishing D}
    Let $P_f$ be a tubular polytope. If $D_{f} \neq 0$, then $\mldeg X_{f} \geq 2$.
\end{lemma}
\begin{proof}
    If $D_{f}$ is not the zero polynomial, then $M_2$ is non-empty and cofinite in particular $\chi(M_2)\leq -1$. As $M_1$ and $M_2$ are disjoint subsets of the torus and $\C^\ast \setminus (M_1 \cup M_2)$ is finite the inequality $\chi(M_1)+ \chi(M_2) \leq 0$ holds.
    It suffices to check the case $\chi (M_2) = -1$. By assumption we know that the polynomial $f_0f_1f_2f_3D_{f}$ has exactly one solution in $\C^\ast$, thus $M_1 = \emptyset$. In particular $\chi (M_1) + \chi (M_2) < 0$. The calculation of the ML degree above yields
    \begin{align*}
        \mldeg X_{f} &= -2\chi(M_2)-\chi(M_1) \\
        & \geq 1 - (\chi(M_2)+\chi(M_1)) \geq 2. \qedhere
    \end{align*}
\end{proof}
\begin{corollary}
   A tubular polytope $P_f$ admits weights with ML degree one if and only if it is square.
\end{corollary}
\begin{proof}
    In the case that $P_f$ is square we define weights by $f_i = x^{a_i}(x-1)^{d_i}$ with $a_i,d_i$ determined by the integer lattice points of the polytope $P_f$, then $M_1 = \C^\ast \setminus \{ 1 \}$ and $M_2 = \emptyset$ give the desired ML degree. On the other hand, if $P_f$ admits weights with ML degree 1, then by the previous lemma $D_f =0$ for those weights. From $f_0f_3 = f_1f_2$ we get trivially that both the order and degree of the polynomials must agree. But those are just the conditions $\ord(f_0)+\ord(f_3) = \ord(f_1)+\ord(f_2)$ and $\deg (f_0) + \deg (f_3) = \deg (f_1) + \deg (f_2)$ hold, i.e., $P_f$ has to be square, we conclude by the previous lemma.
\end{proof}

We would like to better understand the ML degrees that can arise between two and the volume of the polytope for tubular polytopes that have a square base. We see that the formula for the Euler characteristic in terms of $M_1$ and $M_2$ depends on certain root counts and joint
behavior of the roots of the polynomials $f_i$. We now revisit 
Example~\ref{ex:simplex-cut} to illustrate how tweaking the roots lets us observe all possible
ML degrees between two and the normalized volume of the polytope.

\begin{example}\label{ex:simplex-cut2}
Consider again the tubular polytope $P_f$ described in Example~\ref{ex:simplex-cut}. 
We have already exhibited scalings
yielding ML degrees two and three. 
Since the normalized volume of $P_f$ is five, it suffices to find values of the coefficients of $f_i$ that realize ML degrees four and five. To get ML degree four, choose $f_0=1+2x$, $f_1=1$, $f_2=1+9x$, $f_3=1+5x$. In this case $D_f= 2x(5x-1)$, so 
$M_1=\{-1/2,-1/9,-1/5,1/5\}$ and $M_2=\C^{\ast}\setminus \{\textrm{four points}\}$, therefore the ML degree is $-2(-4)-4=4$. Note that $D_f$ is of degree two and it has a zero root that does not count for the Euler characteristic because we are in $\C^{\ast}$.
Changing the choice of polynomials to 
$f_0=1+2x$, $f_1=2$, $f_2=1+5x$, $f_3=1+3x$ and doing the analogous calculations, leads to ML degree 5.
\end{example}

From now on we will assume $\{ w_{000},w_{100},w_{010},w_{110} \}$ to be a facet of $P_f$, i.e. restrict to the case of a tubular polytope with a square base.   
Define $d_A \coloneqq   \deg f_0 + \deg f_3 - \deg f_1 - \deg f_2$. By symmetry, we can assume that $d_A > 0$. Also, since $\{ w_{000},w_{100},w_{010},w_{110} \}$ is a facet of $P_f$, we have $d_A=(d_0+d_3)-(d_1+d_2)=(b_0+b_3)-(b_1+b_2)$.
With this notation, we may write $D_f=x^{a_0+a_3}(g_0g_3 - g_1g_2)$ and $f=x^{a_0}g_0+x^{a_1}g_1s+x^{a_2}g_2t+x^{a_3}g_3st$. Write the polynomials $g_i$ in terms of their factorization $g_i=z_i(1-\frac{x}{\lambda_{1,i}})\cdots (1-\frac{x}{\lambda_{d_i,i}})$ for some tuple $(z_i, \lambda_{1,i},\ldots,\lambda_{d_i,i}) \in (\C^{\ast})^{d_i+1}$. 
Then the scalings of the lattice points in $P_f$ will be given by Laurent polynomials in $z_i, \lambda_{j,i}$. 

Our next goal is to give coefficients for the $f_i$, that achieve each possible ML degree of  the polytope $P_f$. 
The polynomials $f_i$ will be given in terms of the datum $(z_i, \lambda_{1,i},\ldots,\lambda_{d_i,i}) \in (\C^{\ast})^{d_i+1}$ as specified above. In order to define weights for $P_f$ the coefficients of $g_i = \sum_{k=0}^{d_i} (-1)^k e_k(\lambda_1^{-1},\dots,\lambda_{d_i}^{-1})x^k$ need to be nonzero, where $e_k$ denotes the k-th elementary symmetric function. However, for our purpose, we can choose all the roots to be real and positive. 
\begin{lemma}\label{lem:technical}
    Let $n\in \Z_{>0}$ and $\{\alpha_1,\dots ,\alpha_l\} ,\{\beta_1,\dots,\beta_k \}$ be two partitions of $n$. For generic choices of $\lambda_1,\dots,\lambda_l,\mu_2,\dots,\mu_k \in \C^\ast$ the polynomial $H\coloneqq \prod _{i=1}^{l}(1-\frac{x}{\lambda_i})^{\alpha_i} - \prod_{j=2}^k (1-\frac{x}{\mu_j})^{\beta_j}$ has $n-1$ distinct nonzero roots.
\end{lemma}
\begin{proof}
    As $\deg H= n$ and $H(0) = 0$ by construction, for generic choices of $\lambda _i, \mu _j$ we need to show that the discriminant of $H$ is not the zero polynomial in $\lambda_i,\mu_j$. It suffices to show that this holds for $H=(1- \frac{x}{\lambda})^n- 1$. But this can be done by a straightforward computation of solutions to $(1-\frac{x}{\lambda})^n-1 = 0$. Those are given by $x=\lambda(1-\zeta)$ for $\zeta$ a $n$-th root of unity different from 1.
\end{proof}
\begin{lemma}
    The normalized volume of a tubular polytope $P_f$ with square base  is $3(d_0+d_3)-d_A$.
\end{lemma}
\begin{proof}
    One can compute the volume explicitly by counting the number of simplices in a triangulation of $P_f$. To state the triangulation of $P_f$, we recall that its vertices are given by
    \begin{center}
    \begin{tabular}{@{}cccc@{}}
    $w_{000}=(0,0,0),$ &
    $w_{100}=(1,0,0),$ &
    $w_{010}=(0,1,0),$ &
    $w_{110}=(1,1,0),$ \\[4pt]
    $w_{001}=(0,0,d_0),$ &
    $w_{101}=(1,0,d_1),$ &
    $w_{011}=(0,1,d_2),$ &
    $w_{111}=(1,1,d_3).$
    \end{tabular}
    \end{center}
    Define the simplices
    \begin{align*}
        T_1 = &\{ \conv(w_{001},w_{111},(1,0,j),(1,0,j-1)) \colon \, 0 < j \leq d_1 \}, \\
        T_2 = &\{ \conv(w_{001},w_{111},(0,1,j),(0,1,j-1)) \colon \, 0<j\leq d_2  \}, \\
        T_3 = &\{ \conv(w_{001},w_{100},(1,1,j),(1,1,j-1)) \colon \, 0<j\leq d_3\} \\
             \cup &\{\conv(w_{001},w_{010},(1,1,j),(1,1,j-1)) \colon \, 0<j\leq d_3 \}, \\
        T_4 = &\{ \conv(w_{100},w_{110},(0,0,j),(0,0,j-1)) \colon \, 0<j\leq d_0\} \\
             \cup &\{ \conv(w_{010},w_{110},(0,0,j),(0,0,j-1)) \colon \, 0<j\leq d_0\}.
    \end{align*}
    One easily checks that their union equals $P_f$, assuming $d_A\geq 0$. Thus, $T_1 \cup T_2 \cup T_3 \cup T_4 $ is a triangulation of $P_f$ and we obtain
    \begin{align*}
    \vol(P_f) &= \# T_1+\# T_2+\# T_3+\# T_4 = d_1 + d_2 + 2d_3 + 2d_0 = (d_0+d_3-d_A) +2d_3+2d_0 \\
    &= 3(d_0+d_3)-d_A. \qedhere
    \end{align*}
\end{proof}
Next, we want to show that for tubular polytopes with a square base, with the conditions we fixed before Lemma~\ref{lem:technical}, there exist scalings for the lattice points that achieve every ML degree between $d_A+1$ and $\vol(P_f)= 3(d_0+d_3)-d_A$, including both endpoints. 
\begin{proposition}\label{scalings for ML degrees}
    Let $P_f$ be a tubular polytope with a square base. Then each ML degree between $d_A+1$ and $\vol (P_f)$ can be achieved.
\end{proposition}
\begin{proof}
    The strategy is to give explicit scalings for the lattice points of $P_f$ for each ML degree. Again, by symmetry, we make the assumption $d_A > 0$ and also $d_0 > 0$.

    Let $(\alpha _1,\dots, \alpha_r)$ be a partition of $d_0+d_3$, pick 
    \[ 
        (\rho _1, \rho_2,\dots, \rho_{d_0+d_3}) = (\underbrace{\lambda_1,\dots,\lambda_1}_{\alpha_1},\dots,\underbrace{\lambda_r,\dots,\lambda_r}_{\alpha_r})
    \]
    with the $\lambda _i \in \R _{>0}$, $1\leq i \leq r$ chosen generically, in particular they are distinct. Set
    \[
    \begin{aligned}
        g_3&=\prod_{i=1}^{d_3}\left(1-\frac{x}{\rho_i}\right),
        &
        g_0&=\prod_{i=d_3+1}^{d_0+d_3}\left(1-\frac{x}{\rho_i}\right),\\
        g_1&=\prod_{i=1}^{d_1}\left(1-\frac{x}{\rho_i}\right),
        &
        g_2&=\prod_{i=d_1+1}^{d_1+d_2}\left(1-\frac{x}{\rho_i}\right).
    \end{aligned}
    \]
    By construction $D_f = x^{a_0+a_3}g_1g_2(\prod_{i=d_1+d_2+1}^{d_0+d_3}(1-\frac{x}{\rho_i})-1) $, if $d_A > 0$ and  $\lambda_i$ is generic, this has $d_{A}-1$ unique nonzero roots not contained in $V(g_0g_1g_2g_3)$. Let $h = \# ((V(g_0) \cap V(g_3))\setminus V(g_1g_2))$, note that $h\in \{0,1\}$, and $s = \# V(g_1g_2)$, thus
    \[
    \# (\C^\ast \setminus M_2) = r+d_A-1, \quad \# M_1 = r-s-h+d_A-1.
    \]
    The expression of the ML degree can then be given as
    \[
        \mldeg X_f = -2 \chi (M_2) - \chi (M_1) = r+s+h+d_A-1.
    \]
    For $0\leq k < d_1+d_2$ the partition $(\underbrace{1,\dots,1}_{k},d_0+d_3-k)$ gives $k+1$ non-unique nonzero roots for the $f_i$ 
    \[
        \mldeg X_f = 2k+1+d_A
    \]
    and $(\underbrace{1,\dots,1}_{k},d_0+d_3-k-1,1)$ introduces a unique root, so
    \[ 
        \mldeg X_f = 2(k+1) + d_A
    \]
    Thus, we have scalings for ML degrees ranging from $d_A+1$ up to $2(d_1+d_2)+d_A=d_0+d_3+(d_1+d_2)$. 
    Now choose four partitions $(\alpha_1,\dots,\alpha_{r_0}), (\beta_1,\dots,\beta_{r_1}),(\gamma_1, \dots, \gamma_{r_2})$ and $(\delta_1,\dots,\delta_{r_3})$ of $d_0,d_1,d_2$ and $d_3$ respectively. 
    \[
\begin{aligned}
    g_0 &= \prod_{i=1}^{r_0} \left( 1 - \frac{x}{\lambda_i}\right)^{\alpha_i}, &
    g_1 &= \prod_{i=1}^{r_1} \left( 1 - \frac{x}{\mu_i}\right)^{\beta_i} ,\\
    g_2 &= \prod_{i=1}^{r_2} \left( 1 - \frac{x}{\nu _i}\right)^{\gamma_i}, &
    g_3 &= \prod_{i=1}^{r_3} \left( 1 - \frac{x}{\xi _i}\right)^{\delta_i}.
\end{aligned}
    \]
    for generic $\lambda_i, \mu _i, \nu _i, \xi _i \in \R _{>0}$. By Lemma~\ref{lem:technical} $D_f$ has $d_0+d_3-1$ distinct nonzero roots and we compute
    \[
        \mldeg X_f = r_0+r_1+r_2+r_3+d_0+d_3-1.
    \]
 This covers ML degrees $d_0+d_3+\# \{ i \colon d_i > 0 \}-1$ up to $2(d_0+d_3)+d_1+d_2-1=\vol (P_f) - 1$. \\
 We already know the scalings for each ML degree up to $d_0+d_1+d_2+d_3$, trivially $d_1+d_2 \geq \# \{ i\colon d_i>0\}-2$, so the upper bound is at least $d_0+d_3+\# \{ i \colon d_i > 0 \}-2$. \\
 The maximal ML degree $\vol (P_f)$ is achieved for generic nonzero weights on the integer lattice points of $P_f$. 
\end{proof}

\begin{theorem}[Mason--Stothers theorem {\cite{stothers1981}, \cite[Chapter~IV, Theorem~7.1]{lang2002}}]
    Let $a(t),b(t),c(t)$ be relatively prime polynomials, not all constant, and $a+b=c$. Then
    \[
        \max \{ \deg (a) , \deg (b), \deg (c)\} \leq \# V(abc) - 1.   
    \]
\end{theorem}

\begin{proposition}
    The minimal ML degree that can be obtained for a tubular polytope with square base is $d_A+1$.
\end{proposition}
\begin{proof}
    As we have already constructed scalings that achieve ML degree $d_A+1$, it remains to show the lower bound $\mldeg X_f \geq d_A+1$ for $d_A > 0$. 
    
    Let $G \coloneqq \gcd  (g_0g_3,g_1g_2)$, then the assumption of the Mason--Stothers theorem is satisfied for $a=g_0g_3/G$, $b=g_1g_2/G$ and $c=(g_0g_3-g_1g_2)/G$. By assumption $\deg (g_0g_3/G) \geq d_A$, with equality holding if and only if $\deg (g_1g_2/G)=0$. In the case where the inequality is strict, Mason--Stothers yields $\# V(g_0g_1g_2g_3D_f) > \deg (g_0g_3/G) > d_A$, in particular $\#( \C^\ast \setminus M_2 )\geq d_A+1$. But we know that $\# (\C^\ast \setminus M_2)$ is a lower bound for the ML degree. \\

    If $\deg (g_1g_2/G) = 0$, then an application of Mason--Stothers yields $\# (\C^\ast \setminus M_2) \geq d_A$. By the definition of a tubular polytope, $f_1f_2$ cannot be a monomial, or equivalently, $\deg (g_1g_2)>0$. By assumption $\deg(g_1g_2)=\deg (G)$ we have $M_1 \subsetneq \C^\ast \setminus M_2$ and $\mldeg X_f > \# (\C^\ast \setminus M_2 )\geq d_A$, finishing the proof.
\end{proof}

We conclude the section by taking another look at the definition of tubular polytopes; in particular, we have not addressed the reason for the condition $d_0+d_3\geq 1$ and $d_1+d_2\geq 1$. The reason for this assumption is that otherwise $d_A+1$ will not be the lower bound for the ML degree.  If $P_f$ is defined by $f_0,f_1,f_2,f_3$ that are nonzero and not all monomials, but either $f_0f_3$ or $f_1f_2$ is a monomial, then it will be called an exceptional tubular polytope. 

\begin{corollary}
    The ML degree of an exceptional tubular polytope with a square base $P_f$ is at least $d_A$.
\end{corollary}
\begin{proof}
    This is a direct application of the Mason--Stothers theorem for the polynomials $g_0g_3$, $g_1g_2$ and $g_0g_3-g_1g_2$.
\end{proof}

\begin{theorem}[{q-binomial theorem \cite[Eq.~(1.87), p.~74]{stanleyEC1}}]
    For any $j \in \mathbb{N} _{>0}$,
    \[
        \prod_{i=0}^{j-1} (1+xq^i) = \sum_{k=0}^j x^k q^{\binom{k}{2}}
            \begin{pmatrix}
                j \\ k
            \end{pmatrix}_{q}
    \]
    with the q-binomial coefficient defined as
    \[
        \begin{pmatrix}
            j \\ k
        \end{pmatrix}_q
        = \frac{(1-q^j)(1-q^{j-1})\hdots(1-q^{j-k+1})}{(1-q)(1-q^2)\hdots(1-q^k)}.
    \]
\end{theorem}

\begin{proposition}
    For an exceptional tubular polytope with square base $P_f$, each ML degree between $d_A$ and $\vol (P_f)$ can be achieved.
\end{proposition}
\begin{proof}
    By defining the $f_i$ of exceptional tubular polytopes with a square base as in Proposition~\ref{scalings for ML degrees} we can cover all cases but one. If $d_1=d_2=0$ but both $d_0$ and $d_3$ are positive, we need a scaling that gives the ML degree $d_A=d_0+d_3$. Let $\zeta$ be a primitive $d_A$-th root of unity
    \[
    \begin{aligned}
        g_0&=\prod_{i=0}^{d_0-1}\left(1-\frac{x}{\zeta ^i}\right),
        &
        g_1&=1,
        &
        g_2&=1,
        &
        g_3&=\prod_{i=0}^{d_3-1}\left(1-\frac{x}{\zeta ^{d_0+i}} \right).
    \end{aligned}
    \]
    Then $D_f = x^{a_0+a_3}(\prod_{i=0}^{d_A-1} (1-\frac{x}{\zeta ^i}) -1)= -x^{a_0+a_3+d_A}$ and $\mldeg X_f = d_A$. Up to sign, the coefficients of $g_0$ are given by the q-binomial theorem evaluated at $q= \zeta$. As $\zeta$ is a primitive root of unity and $d_0<d_A$ none of the q-binomial coefficients vanishes. \\
    A substitution of $x$ by $\zeta^{d_0}y$ gives $g_3= \prod_{i=0}^{d_3-1}(1-\frac{y}{\zeta ^i})$ we can conclude with the same argument.
\end{proof}
We end this section with one example, using the polytope in Example~\ref{ex:simplex-cut}, explaining how to  construct a  closed set $S_2$ such that every scaling $c\in S_2$ restricts to ML degree one on the faces of the polytope.

\begin{example} \label{ex:cube-simplex}
The vertices of the polytope in Example~\ref{ex:simplex-cut}
are the columns of the matrix
\[
    A=\begin{pmatrix}
         0& 1& 0& 1& 0& 0& 1\\
         0& 0& 1& 1& 0& 1& 1\\
         0& 0& 0& 0& 1& 1& 1\\
    \end{pmatrix}
    \]
     We associate the vector of scalings $c=(c_1,c_2,c_3,c_4,c_5,c_6,c_7)$ to the columns of $A$.
     The polynomial associated to $(A,c)$ and its homogenization are \begin{align*}
        f&=c_1+c_2s +c_3t+c_4st+c_5u+c_6tu+c_7stu \\
        F&=c_1W^3+c_2XW^2 +c_3YW^2+c_4XYW+c_5ZW^2+c_6YZW+c_7XYZ.
    \end{align*}
We describe the locus of points $c\in(\mathbb{C}^{\ast})^7$ such that for each facet $H$ of $P_f$, the restriction
    $(H,c_{H})$ has ML degree one. The polytope $P_f$ has three square faces and the rest are triangles, these triangles
    have ML degree one for any choice of scaling, thus we focus only on the square faces. To identify the necessary conditions on $c$ to preserve the ML degree one property on the faces, number the faces by
    $
        F_1=\conv(a_1,a_2,a_3,a_4),
        F_2=\conv(a_3,a_4,a_6,a_7),
        F_3= \conv(a_1,a_3,a_5,a_6),
  $
    where the $a_i$ indicates the column $i$ of the matrix $A$. The restriction of the polynomial $F$ to the face $F_i$
    is denoted $G_i$. We obtain three polynomials $G_1,G_2,G_3$ and
    from \cite{vonBothmer2010}, we know that the degree two factor of each of the polynomials $G_i$ has ML degree one if and only if it factors into a product of two linear forms. The choice of coefficients in this
    factorization is not important and we may assume it is as follows:
    \begin{align*}
        G_1&=W(c_1W^2+c_2XW+c_3YW+c_4XY)=W(a_1W+a_2X)(a_3W+a_4Y)\\
        G_2&=Y(c_3W^2+c_4XW+c_6ZW+c_7XZ)=Y(b_1W+b_2X)(b_3W+b_4Z) \\
        G_3&=W(c_1W^2+c_3YW+c_5ZW+c_6YZ)=W(d_1W+d_2Y)(d_3W+d_4Z)
    \end{align*}
for some nonzero $a_i,b_i,d_i$. Expanding out the product and setting coefficients of monomials in $X,Y,Z$ equal to each
other, yields an ideal generated by the following polynomials.
\begin{align*}
    c_1-&a_1a_3,\quad c_3-a_1a_4,\quad c_2-a_2a_3,\quad c_4-a_2a_4, \quad c_3-b_1b_3,\quad c_6-b_1b_4, \\ 
c_4-&b_2b_3,\quad c_7-b_2b_4,\quad c_1-d_1d_3,\quad c_5-d_1d_4,\quad c_3-d_2d_3,\quad c_6-d_2d_4.
\end{align*}
Eliminating the variables  $a_i,b_i,d_i$, and saturating by the
$c_i$ yields an ideal defined by
\[I=\langle c_4c_6 - c_3c_7, c_2c_6 - c_1c_7, c_4c_5 - c_1c_7, c_3c_5 - c_1c_6, c_2c_3 - c_1c_4 \rangle.
\]
Any point in $S_2=V(I)\subset (\C^{\ast})^7$, when restricted to any one of the facets, has ML degree one. Recall that in Example~\ref{ex:simplex-cut} and
Example~\ref{ex:simplex-cut2} we constructed scalings that give ML degree 2,3,4,5. We can substitute these scalings in the ideal to understand which ones restrict to ML degree one on the faces. We summarize this in a table.
\[
\begin{array}{c|ccccccc|c}
\text{ML degree} & c_1 & c_2 & c_3 & c_4 & c_5 & c_6 & c_7 & \text{ML degree one on faces} \\
\hline
2 & 1 & 1 & 1 & 1 & 1 & 1 & 1 & \text{yes} \\
3 & 1 & 2 & 1 & 2 & 3 & 3 & 2 & \text{no}  \\
4 & 1 & 1 & 1 & 1 & 2 & 9 & 5 & \text{no}  \\
5 & 1 & 2 & 1 & 1 & 2 & 5 & 3 & \text{no}
\end{array}
\]
\end{example}

\section*{Acknowledgements}
The authors thank the Simons Laufer Mathematical Sciences Institute (SLMath) and the
Max Planck Institute for Mathematics in the Sciences for hosting the summer school
\emph{New Perspectives on Discriminants and Their Applications}, where this work began.
Eliana Duarte was partially supported by CMUP, member of LASI, which is financed by
Portuguese national funds through FCT -- Fundação para a Ciência e a Tecnologia, I.P.,
under the project with reference UID/00144/2025 and DOI
\url{https://doi.org/10.54499/UID/00144/2025}.
\section*{Declaration of generative AI use}
During the final preparation of this manuscript, the authors used Claude Opus 5.5 and GPT-Astra 6 to assist with proofreading and in preparing Figure~\ref{fig:Pf}. The models also pointed us to the
Mason--Stothers theorem and the $q$-binomial theorem. All mathematical content, including statements and
proofs, was written and checked by the authors, who take full responsibility for
the contents of this paper.
\bibliographystyle{amsplain}
\bibliography{refs}
\end{document}